\documentclass[11pt]{article}
\usepackage{amsfonts}
\usepackage{latexsym, amssymb, amsmath, graphicx}
\usepackage{graphicx}

\newtheorem{thm}{Theorem}[section]

\newtheorem{defi}[thm]{Definition}
\newtheorem{lem}[thm]{Lemma}

\def\pf{\noindent{\it Proof.} }
\def\qed{\nopagebreak\hfill{\rule{4pt}{7pt}}\medbreak}

\numberwithin{equation}{section}

\def\qed{\nopagebreak\hfill{\rule{4pt}{7pt}}
\medbreak}

\title{The Rogers-Ramanujan-Gordon Theorem for Overpartitions }
\author{William Y.C. Chen\raisebox{5pt}{\scriptsize 1},
Doris D. M. Sang\raisebox{5pt}{\scriptsize 2}, and Diane Y. H.
Shi\raisebox{5pt}{\scriptsize 3}}
\date{Center for Combinatorics, LPMC-TJKLC\\
 Nankai University\\
Tianjin 300071, P.R. China \\
\vspace{15pt} \raisebox{5pt}{\scriptsize 1\,}chen@nankai.edu.cn,
\raisebox{5pt}{\scriptsize 2\,}sdm@cfc.nankai.edu.cn,
\raisebox{5pt}{\scriptsize 3\,}shiyahui@cfc.nankai.edu.cn}

\begin{document}

\allowdisplaybreaks

%--------------------------------------------------------------------
\maketitle
%--------------------------------------------------------------------
\noindent {\bf Abstract.}
Let $B_{k,i}(n)$ be the number  of partitions of $n$ with
certain difference condition and let $A_{k,i}(n)$ be the number of
partitions of $n$ with certain congruence condition. The Rogers-Ramanujan-Gordon theorem  states that $B_{k,i}(n)=A_{k,i}(n)$.
   Lovejoy obtained an overpartition analogue of the Rogers-Ramanujan-Gordon theorem
   for the cases $i=1$ and $i=k$. We find an overpartition analogue of
   the Rogers-Ramanujan-Gordon theorem in the general case. Let $D_{k,i}(n)$ be the number  of overpartitions of $n$ satisfying  certain difference condition and $C_{k,i}(n)$ be the number of overpartitions of $n$  whose non-overlined parts satisfy certain congruences condition. We show that $C_{k,i}(n)=D_{k,i}(n)$.
   By using a function introduced by Andrews,
   we obtain a recurrence relation which implies that
   the generating function of $D_{k,i}(n)$ equals the generating function of $C_{k,i}(n)$.
   We also find a generating function formula
    of $D_{k,i}(n)$ by using   Gordon marking representations of overpartitions,
    which can be considered as an overpartition analogue of an identity of Andrews
    for ordinary partitions.


%--------------------------------------------------------------------

\section{ Introduction}

In this paper, we obtain the  Rogers-Ramanujan-Gordon theorem
 for overpartitions. Furthermore, by introducing the Gordon marking of
  an overpartition, we find a generating function
formula which can be considered as an overpartition analogue of an identity of Andrews.
Notice that the identity of Andrews implies the Rogers-Ramanujan-Gordon theorem
for ordinary partitions, see Kur\c{s}ung\"{o}z \cite{kur09}.

  An overpartition is a partition for which the first occurrence of a part may be overlined. For example, $(\overline{7},7,6,\overline{5},2,\overline{1})$ is an overpartition of $28$.  There are many $q$-series identities that have combinatorial
interpretations in terms of overpartitions,
see, for example, Corteel and Lovejoy \cite{cor04}. Furthermore, overpartitions possess many analogous properties of ordinary partitions, see Lovejoy \cite{lov03, lov06}. For example,
various   overpartition analogues of the Rogers-Ramanujan-Gordon theorem have been
obtained by  Corteel and
Lovejoy \cite{cor07}, Corteel, Lovejoy and Mallet \cite{cor08} and Lovejoy \cite{lov03, lov04, lov07, lov10}.

Let us recall that Gordon \cite{gor61} found the following combinatorial generalization of the
Rogers-Ramanujan identities \cite{rog1894}, which has been called the Rogers-Ramanujan-Gordon
theorem, see Andrews \cite{and66}.

\begin{thm}(Rogers-Ramanujan-Gordon)\label{thm1.2} Let $B_{k,i}(n)$ denote the number of partitions
of $n$ for the form $b_1 + b_2 + \cdots + b_s$, where $b_j \geq b_{j+1}$,
$b_j-b_{j+k-1}\geq2$ and at most $i-1$ of the $b_j$ are equal to $1$
and $1\leq i\leq k$. Let $A_{k,i}(n)$ denote the number of
partitions of $n$ into parts $\not \equiv0,\pm i (mod\ 2k + 1)$.
Then for all $n \geq0$, \[A_{k,i}(n) = B_{k,i}(n).\]
\end{thm}

Lovejoy \cite{lov03} obtained overpartition analogues of the above
Rogers-Ramanujan-Gordon  theorem for $i=k$ and $i=1$.

\begin{thm}\label{lov1}
Let $\overline{B}_{k}(n)$ denote the number of
overpartitions of $n$ of the form $y_1+y_2+\cdots+y_s$, such that
$y_j-y_{j+k-1}\geq1$ if $y_j$ is overlined and $y_j-y_{j+k-1}\geq2$
otherwise. Let $\overline{A}_{k}(n)$ denote the number of overpartitions of
$n$ into parts not divisible by $k$. Then $\overline{A}_{k}(n)=\overline{B}_{k}(n)$.
\end{thm}

\begin{thm}\label{lov2} Let $\overline{D}_{k}(n)$ denote the number of overpartitions of $n$ of the
form $z_1+z_2+\cdots+z_s$, such that $1$ cannot occur as a
non-overlined part, and where $z_j-z_{j+k-1}\geq1$ if $z_j$ is
overlined and $z_j-z_{j+k-1}\geq2$ otherwise. Let $\overline{C}_{k}(n)$ denote
the number of overpartitions of $n$ whose non-overlined parts are
not congruent to $0,\pm1$ modulo $2k$. Then $\overline{C}_{k}(n)=\overline{D}_{k}(n)$.
\end{thm}

The first  result of this paper
is to give an overpartition analogue of
the Rogers-Ramanujan-Gordon theorem in the general case.

\begin{thm}\label{thm1} For $k\geq i\geq 1$, let $D_{k,i}(n)$ denote the number of overpartitions of $n$ of the form $d_1+d_2+\cdots+d_s$,
such that $1$ can occur as a non-overlined part at most $i-1$ times,
and where $d_j-d_{j+k-1}\geq1$ if $d_j$ is overlined and
$d_j-d_{j+k-1}\geq2$ otherwise. For $k> i\geq 1$, let $C_{k,i}(n)$ denote the number of
overpartitions of $n$ whose non-overlined parts are not congruent to
$0,\pm i$ modulo $2k$ and let $C_{k,k}(n)$ denote the number of overpartitions of $n$ with parts not divisible by $k$.
 Then $C_{k,i}(n)=D_{k,i}(n)$.
 \end{thm}

It is clear that  Theorem \ref{thm1} contains Theorems \ref{lov1} and \ref{lov2}  as special
cases for $i=k$ and $i=1$. To be more specific,
$\overline{B}_{k}(n)$ and $\overline{A}_k(n)$ in Theorem \ref{lov1} are $D_{k,k}(n)$ and $C_{k,k}(n)$ in Theorem \ref{thm1},  $\overline{D}_{k}(n)$ and $\overline{C}_k(n)$ in Theorem \ref{lov2} are
$D_{k,1}(n)$ and $C_{k,1}(n)$ in Theorem \ref{thm1}.

We will give an algebraic proof of  Theorem \ref{thm1} in the next section by
showing that the generating function of $D_{k,i}(n)$ equals the generating function of  $C_{k,i}(n)$. It is evident that the generating function of $C_{k,i}(n)$ equals
 \begin{equation}\label{cki}\sum_{n\geq 0}C_{k,i}(n)q^n=\frac{(-q)_\infty(q^i,q^{2k-i},q^{2k};q^{2k})_\infty}{(q)_\infty}.
\end{equation}
In fact, we shall prove a stronger result on a refinement  of the generating function
of $D_{k,i}(n)$.

The generating function versions  of Theorem \ref{thm1.2} for $k=2$ are the  Rogers-Ramanujan identities
\begin{equation}\label{rr1}
\sum_{n\geq0}\frac{q^{n^2+n}}{(q)_n}=\frac{1}{(q^2,q^3;q^5)_\infty},
\end{equation}
and
\begin{equation}\label{rr2}
\sum_{n\geq0}\frac{q^{n^2}}{(q)_n}=\frac{1}{(q^1,q^4;q^5)_\infty}.
\end{equation}

Note  that the left hand sides of \eqref{rr1} and \eqref{rr2} can be  interpreted as the
generating functions for $B_{2,1}(n)$ and $B_{2,2}(n)$ respectively.
As a generalization of the Rogers-Ramanujan identities,
 Andrews \cite{and74} obtained the following theorem.
 \begin{thm}\label{and}For $k\geq i\geq1$,
 \begin{equation}\label{equ2}
\sum_{N_1\geq N_2\geq\cdots\geq
N_{k-1}\geq0}
\frac{q^{N_1^2+N_2^2+\cdots+N_{k-1}^2+N_{i}+\cdots+N_{k-1}}
}{(q)_{N_1-N_2}\cdots(q)_{N_{k-2}-N_{k-1}}(q)_{N_{k-1}}}=
\frac{(q^i,q^{2k+1-i},q^{2k+1};q^{2k+1})_\infty}{(q)_\infty}.
\end{equation}
\end{thm}

The sum on the left hand side of (\ref{equ2})
can be viewed as the generating function
for $B_{k,i}(n)$.
Andrews  proved that the both sides
of \eqref{equ2} satisfy the same recurrence relation.

While it is easy to give
combinatorial interpretations of the left hand sides of \eqref{rr1} and \eqref{rr2},
it does not seem to be trivial to show that
the left hand side of \eqref{equ2} is  the generating function for $B_{k,i}(n)$.
Kur\c{s}ung\"{o}z \cite{kur09} provided a combinatorial explanation
of the left hand side of \eqref{equ2}  by introducing the notion of the Gordon marking of a
partition. More precisely, he obtained the  following formula for the generating function
of $B_{k,i}(m,n)$, where $B_{k,i}(m,n)$ denotes the number of partitions enumerated by $B_{k,i}(n)$
that have $m$ parts.

\begin{thm}\label{Bki}For $k\geq i\geq 1$, \begin{align}\label{Bmn}
\sum_{m,n\geq 0}B_{k,i}(m,n)x^mq^n
=\sum_{N_1\geq\cdots\geq
N_{k-1}\geq0}
\frac{q^{N_1^2+N_2^2+\cdots+N_{k-1}^2+N_{i}+\cdots+N_{k-1}}
x^{N_1+\cdots+N_{k-1}}}{(q)_{N_1-N_2}\cdots(q)_{N_{k-2}-N_{k-1}}(q)_{N_{k-1}}}.
\end{align}
\end{thm}

The second result of this paper is the following formula for the
generating function of the number $D_{k,i}(m,n)$  of overpartitions enumerated by $D_{k,i}(n)$
that have $m$ parts. We shall give a combinatorial proof of this identity by using the
Gordon marking representations of overpartitions.

\begin{thm}\label{main} For $k\geq i\geq 1$, we have
\begin{align}\label{equ1}
&\sum_{n=0}^{\infty}D_{k,i}(m,n)x^mq^n\nonumber
\\& \quad
=\sum_{N_1\geq \cdots\geq
N_{k-1}\geq0}\frac{q^{\frac{(N_1+1)N_1}{2}+N_2^2+\cdots+N_{k-1}^2+N_{i+1}+\cdots+N_{k-1}}
(-q)_{N_1-1}(1+q^{N_i})x^{N_1+\cdots+N_{k-1}}}{(q)_{N_1-N_2}\cdots(q)_{N_{k-2}-N_{k-1}}(q)_{N_{k-1}}},\end{align}
where assume that $N_k=0$.
\end{thm}

By setting $x=1$ in \eqref{equ1}, we obtain the generating function for $D_{k,i}(n)$ which is the left hand side
of \eqref{a}.
By  Theorem \ref{thm1}, we are led to the following theorem which can be seen as an overpartition analogue of Andrews' identity
\eqref{equ2}.
\begin{thm}\label{cd}For $k\geq i\geq1$,
\begin{align}\label{a}
&\sum_{N_1\geq\cdots\geq N_{k-1}\geq0}
\frac{q^{\frac{(N_1+1)N_1}{2}+N_2^2+\cdots+N_{k-1}^2+N_{i+1}+\cdots+N_{k-1}}
(-q)_{N_1-1}(1+q^{N_i})}{(q)_{N_1-N_2}\cdots(q)_{N_{k-2}-N_{k-1}}(q)_{N_{k-1}}}\nonumber
\\[6pt] &  \qquad \qquad =\frac{(-q)_\infty(q^i,q^{2k-i},q^{2k};q^{2k})_\infty}{(q)_\infty}.
\end{align}
\end{thm}

It is clear that the generating function for $C_{k,i}(n)$ equals the right hand  side of \eqref{a}. Hence identity \eqref{a}  can be viewed as the generating function version of Theorem \ref{thm1}.
It should be noticed that the approach of Andrews to
 \eqref{equ2} for ordinary partitions does not seem to apply to the
 above identity \eqref{a} for overpartitions.

The special case of  identity \eqref{a} for $i=1$ was derived by Chen, Sang and Shi
\cite{chen} by using Andrews' multiple series transformation  \cite{and75}. In this case, the left hand side of \eqref{a} has a combinatorial interpretation in terms of the generating function of  the number of anti-lecture hall compositions of $n$
with the first entry not exceeding $2k-2$.

The special case  of \eqref{a}
for $i=k$ was obtained by Corteel and Lovejoy \cite{cor04} also by using Andrews' multiple series transformation. In this case, the left hand side of \eqref{a}
has a combinatorial interpretation  in terms of
 the number of overpartitions
whose Frobenius representation has a top row with at most $k-2$ Durfee squares in its associated partition.

 However, for $2\leq i\leq k-1$,  identity (\ref{a}) does not seem to be
a consequence of Andrews' multiple series transformation.
It should be mentioned that for  $i=1,k$,
 the combinatorial interpretation of the left hand side of
(\ref{a}) as a Rogers-Ramanujan-Gordon theorem for overpartitions  as in Theorem \ref{thm1}
 is different from the interpretation in terms
of anti-lecture hall compositions given in Chen, Sang and Shi \cite{chen} or the Frobenius representations given in Corteel and Lovejoy \cite{cor04}.

This paper is organized as follows. In Section 2, we  give an algebraic proof of Theorem \ref{thm1} by showing that
 $C_{k,i}(n)$ and $D_{k,i}(n)$
 satisfy the same recurrence relation.
In Section 3, we introduce the notion of the
Gordon marking of an  overpartition.
To  prove Theorem \ref{main}, we divide the set of overpartitions enumerated
 by $D_{k,i}(m,n)$ into two subsets.
In Section 4, we define the first reduction operation and the first dilation operation. Based on the these two operations we give the first bijection for the proof of Theorem \ref{main}.  In Section 5, we introduce the second reduction operation and the second dilation operation on the Gordon marking representations of
overpartitions. Then we give the second bijection for the proof of Theorem \ref{main}.
In Section 6, we give the third  bijection for the proof of Theorem \ref{main}.  In Section 7, we complete the proof of the Theorem \ref{main}.

\section{An algebraic proof of Theorem \ref{thm1}}

In this section, we give an algebraic proof of Theorem \ref{thm1},
that is, $C_{k,i}(n)=D_{k,i}(n)$ for any $k\geq i\geq 1$.
We shall use a series $H_{k,i}(a;x;q)$  introduced by
 Andrews \cite{and66, and74}, which is defined by
\begin{equation}\label{eqH}
H_{k,i}(a;x;q)=
\sum_{n=0}^{\infty}\frac{x^{kn}q^{kn^2+n-in}a^n(1-x^iq^{2ni})
(axq^{n+1})_{\infty}(1/a)_n}{(q)_n(xq^n)_\infty}.
\end{equation}
In his algebraic proof of the Rogers-Ramanujan-Gordon theorem,
Andrews used the function $J_{k,i}(a;x;q)$ constructed based on $H_{k,i}(a;x;q)$,
\begin{equation} \label{jki}
J_{k,i}(a;x;q)=
H_{k,i}(a;xq;q)-axqH_{k,i-1}(a;xq;q).
\end{equation}

Lovejoy \cite{lov06} proved Theorem \ref{lov1} and Theorem \ref{lov2} also by using $J_{k,i}(a;x;q)$
for special values of $a$ and $x$.
More precisely, he showed
 the generating function of $\overline{A}_{k}(n)$ and $\overline{C}_{k}(n)$, namely, $C_{k,k}(n)$ and $C_{k,1}(n)$,
 are given by the functions
 $J_{k,k}(-1;1;q)$ and $J_{k,1}(-1/q;1;q)$.
As pointed out by Lovejoy,
the approach of using the function $J_{k,i}(a;x;q)$ does not seem to apply to the general
case, since for $i\not=1, k$, the functions $J_{k,i}(-1;1;q)$ and $J_{k,i}(-1/q;1;q)$ do not
appear to be expressible as single infinite products.

We find that for overpartitions
the function $H_{k,i}(a;x;q)$ itself is the right choice
to prove that $C_{k,i}(n)=D_{k,i}(n)$ for all $k\geq i\geq1$.
 In fact, we shall show that the generating function
 of $C_{k,i}(n)$ can be expressed in terms of $H_{k,i}(a;x;q)$ for special values of
 $a$ and $x$. To explain the fact that the generating functions of $C_{k,k}(n)$ and $C_{k,1}(n)$
 can also be expressed by $J_{k,k}(-1;1;q)$ and $J_{k,1}(-1/q;1;q)$, we have the observations
 \begin{equation}
J_{k,k}(-1;1;q)=H_{k,k}(-1/q;q;q),
 \end{equation}
 and
 \begin{equation}
 J_{k,1}(a;x;q)=H_{k,1}(a;xq;q).
 \end{equation}

Andrews \cite{and66,and76} showed that the generating function  of $B_{k,i}(m,n)$  can be expressed by $J_{k,i}(a;x;q)$:
\begin{equation}\sum_{m,n\geq 0}B_{k,i}(m,n)x^mq^n=J_{k,i}(0;x;q).
\end{equation}
We shall give the following theorem which involves  a refinement of  the number
$D_{k,i}(n)$. Recall that $D_{k,i}(m,n)$ is the number of overpartitions enumerated by $D_{k,i}(n)$ with $m$ parts. As will be seen, once the generating function of $D_{k,i}(m,n)$
is obtained, it is easy to derive the generating function of $D_{k,i}(n)$ by using
Jacobi's triple product identity.

\begin{thm}\label{DH}For $k\geq i\geq 1$, we have
\begin{equation}\sum_{m,n\geq 0}D_{k,i}(m,n)x^mq^n=H_{k,i}(-1/q;xq;q).
\end{equation}
\end{thm}

\noindent {\it Proof.}
  We define
  \begin{equation}\label{wkid}
   W_{k,i}(x;q)=H_{k,i}(-1/q;xq;q),
   \end{equation}
  and \begin{equation}
  \label{wki}W_{k,i}(x;q)=\sum_{m,n=-\infty}^{\infty}W_{k,i}(m,n)x^mq^n.
  \end{equation}
  By the recurrence relation of $H_{k,i}(a;x;q)$,  one can derive
   a recurrence relation of $W_{k,i}(m,n)$. It is easy to give a combinatorial interpretation of $D_{k,i}(m,n)-D_{k,i-1}(m,n)$. This yields  a recurrence relation of $D_{k,i}(m,n)$ which
    coincides with a recurrence relation of $W_{k,i}(m,n)$.

Recall that $H_{k,i}(a;x;q)$ satisfies the following recurrence relation, see Andrews \cite[Lemma 7.1]{and76},
\begin{equation}\label{hki}H_{k,i}(a;x;q)-H_{k,i-1}(a;x;q)
=x^{i-1}H_{k,k-i+1}(a;xq;q)-ax^iqH_{k,k-i}(a;xq;q).
\end{equation}
Substituting  $a=-1/q$ and $x=xq$ into (\ref{hki}), we obtain
\begin{equation}\label{eqJ}
W_{k,i}(x;q)-W_{k,i-1}(x;q)=(xq)^iW_{k,k-i}(xq;q)+(xq)^{i-1}W_{k,k-i+1}(xq;q).
\end{equation}

Our goal is to prove
that  $D_{k,i}(m,n)$  equals $W_{k,i}(m,n)$.
In doing so, we shall show that $D_{k,i}(m,n)$ and $W_{k,i}(m,n)$ satisfy the
same recurrence relation with the same initial values, where $W_{k,i}(m,n)$
is the coefficient of $x^mq^{n}$ in the expansion of $W_{k,i}(x;q)$,
as given by \eqref{wki}.

 Clearly, we have the initial values
 $W_{k,i}(0,0)=1$ for $ k\geq i\geq 1$ and
 $W_{k,0}(m,n)=0$ for $k\geq 1, m,n\geq 0$. Moreover, we assume that
 $W_{k,i}(m,n)=0$ if $m$ or $n$ is zero but not both, and $W_{k,i}(m,n)=0$
 if $m$ or $n$ is negative.
From \eqref{eqJ} it is easily seen that
\begin{equation}
\label{j4}W_{k,i}(m,n)-W_{k,i-1}(m,n)=W_{k,k-i}(m-i,n-m)+W_{k,k-i+1}(m-i+1,n-m),
\end{equation}
Thus  $W_{k,i}(m,n)$ can be defined by the recurrence relation
\eqref{j4} along with the initial values.

Next we wish to find a recurrence relation of $D_{k,i}(m,n)$.
It can be verified that $D_{k,i}(m,n)$ has the initial values
 $D_{k,i}(0,0)=1$ for $ k\geq i\geq 1$ and
 $D_{k,0}(m,n)=0$ for $k\geq 1, m,n\geq 0$. Clearly,  if exactly one of  $m$ and $n$
  is zero, then $D_{k,i}(m,n)=0$. If one of $m$ and $n$ is negative, then
  $D_{k,i}(m,n)=0$. Hence $D_{k,i}(m,n)$ has the same initial values as $W_{k,i}(m,n)$.
It remains to prove that
\begin{equation}\label{dkir}
D_{k,i}(m,n)-D_{k,i-1}(m,n)=D_{k,k-i}(m-i,n-m)+D_{k,k-i+1}(m-i+1,n-m).
\end{equation}

From the definition of $D_{k,i}(m,n)$, one sees that
$D_{k,i}(m,n)-D_{k,i-1}(m,n)$ equals the number of overpartitions
enumerated by $D_{k,i}(m,n)$ such that the non-overlined part $1$
appears exactly $i-1$ times.  We shall divide  the overpartitions enumerated by
$D_{k,i}(m,n)-D_{k,i-1}(m,n)$ into two classes so that we can give
a combinatorial interpretation of the right hand side of (\ref{dkir}).

Let $S_1$ be the set of overpartitions
  enumerated by $D_{k,i}(m,n)-D_{k,i-1}(m,n)$ that contain a part $\overline{1}$,
   and let $S_2$ be the set of overpartitions enumerated by $D_{k,i}(m,n)-D_{k,i-1}(m,n)$
   that do not contain the part $\overline{1}$.
We shall show that the number of  overpartitions in $S_1$  equals $D_{k,k-i}(m-i,n-m)$ and
the number of overpartitions in $S_2$ equals $D_{k,k-i+1}(m-i+1,n-m)$.

Let $\lambda$ be an overpartition in $S_1$.
So $\lambda$ has $i$ parts equal to $1$ or $\overline{1}$.
Removing these $i$ parts, we obtain an overpartition that   contains neither $1$ nor $\overline{1}$.
Subtracting $1$ from each part of the resulting overpartition,
we get an overpartition $\lambda'$. More precisely, by subtracting $1$ from $\overline{r}$
we mean to change $\overline{r}$ to $\overline{r-1}$. From the definition of $D_{k,i}(m,n)$, we find that the parts $1$, $\overline{1}$ and $2$ occur at most $k-1$ times. Notice
that the number of occurrences of
$1$ and $\overline{1}$ in $\lambda$ equals $i$.  Thus,  $2$
appear at most $k-i-1$ times in $\lambda$. So after the
subtraction, the part $1$ appears at most $k-i-1$ times in $\lambda'$.
 By the definition of $D_{k,i}(m,n)$, we deduce that the resulting overpartition $\lambda'$ is
enumerated by $D_{k,k-i}(m-i,n-m)$. Moreover, it is readily seen  that every overpartition enumerated by $D_{k,k-i}(m-i,n-m)$ can be constructed by the above procedure.

For an overpartition  $\lambda$ in $S_2$,
there are exactly $i-1$ parts equal to $1$ in $\lambda$, so
the part  $2$  occurs at most $k-i$ times in $\lambda$.
Removing the $i-1$ parts $1$ and subtracting $1$ from each of
the remaining parts,  we get an overpartition $\lambda'$.
It can be seen that the part $1$ appears $k-i$ times in  $\lambda'$. By the definition of $D_{k,k-i+1}(m-i+1,n-m)$, we find that $\lambda'$ is enumerated by $D_{k,k-i+1}(m-i+1,n-m)$.
Conversely,  every overpartition enumerated by $D_{k,k-i+1}(m-i+1,n-m)$ can be constructed from an overpartition $\lambda$ in $S_2$.

So we have proved  relation \eqref{dkir}, which implies that
$D_{k,i}(m,n)=W_{k,i}(m,n)$ for all $k\geq i\geq 1$, and $m,n\geq 0$,
since $D_{k,i}(m,n)$ and $W_{k,i}(m,n)$ have the same initial values. Thus the generating function of $D_{k,i}(m,n)$  equals $W_{k,i}(x;q)$.  This completes the proof. \qed

We are ready to prove Theorem \ref{thm1}.
Let us compute  the generating function  of $D_{k,i}(n)$.
Setting $x=1$ in Theorem \ref{DH}, we obtain that
\begin{align*}
H_{k,i}(-1/q;q;q) & =\sum_{n=0}^{\infty}\frac{(-1)^nq^{kn^2+kn-in}
(1-q^{(2n+1)i})(-q^{n+1})_{\infty}(-q)_n}
{(q)_n(q^{n+1})_{\infty}}\\[5pt]
&=\frac{(-q)_{\infty}}{(q)_{\infty}}\sum_{n=0}^{\infty}(-1)^nq^{kn^2+kn-in}(1-q^{(2n+1)i})
\\[5pt]
&=\frac{(-q)_{\infty}}{(q)_{\infty}}\sum_{n=-\infty}^{\infty}(-1)^nq^{kn^2+kn-in}.
\end{align*}
In view of Jacobi's triple product identity, we find that
\begin{equation}
\sum_{n\geq 0}D_{k,i}(n)q^n
=\frac{(q^i,q^{2k-i},q^{2k};q^{2k})_{\infty}(-q)_{\infty}}{(q)_{\infty}},
\end{equation}
which implies that $C_{k,i}(n)=D_{k,i}(n)$. This completes the proof of Theorem \ref{thm1}.

\section{The Gordon marking of an overpartition}

 In this section, we introduce the notion of the
 Gordon marking of an overpartition and  give an outline of the proof
 of the generating function formula for $D_{k,i}(m,n)$ as   stated in Theorem \ref{main}.
 To compute the generating function of
 $D_{k,i}(m,n)$, we  divide the set enumerated by $D_{k,i}(m,n)$ into two classes $U_{k,i}(m,n)$
 and $I_{k,i}(m,n)$. Let $F_{k,i}(m,n)$ be the
 number of overpartitions in $U_{k,i}(m,n)$.
 By two simple bijections we can express the generating function of $D_{k,i}(m,n)$ by the generating function of
$F_{k,i}(m,n)$. We shall give the generating function of $F_{k,i}(m,n)$ in Theorem \ref{F}. As will be seen,
we need three bijections to prove Theorem \ref{F},
 which will be presented in Sections 4--6.

 Notice that  identity \eqref{equ2} of Andrews \cite{and74} is a generalization of the Rogers-Ramanujan identity.
 It is natural to ask whether there is  an overpartition analogue
 of  \eqref{equ2}. The answer is given in  Theorem \ref{cd}.
To this end, we  shall give a combinatorial treatment of the generating function of $D_{k,i}(m,n)$
by introducing the notion of Gordon marking representations of overpartitions.
 Observe
that the generating function of $D_{k,i}(m,n)$ stated in Theorem \ref{main} is in the form of  the left hand
side of  \eqref{equ2}.  Thus Theorem \ref{cd} can be deduced from
 Theorem \ref{main} and Theorem \ref{thm1}.

 Kur\c{s}ung\"{o}z \cite{kur09} introduced the notion of
 the Gordon marking of an ordinary
 partition and gave a combinatorial interpretation of identity \eqref{Bmn}.
A Gordon marking of an ordinary partition $\lambda$ is an assignment of
positive integers (marks) to parts of $\lambda$ such that any two equal parts,
as well as
any two nearly equal parts $j$ and $j+1$ are assigned different marks,
 and the marks are as small as possible assuming that the
marks are assigned to the parts in increasing order. For example, the Gordon marking of
\[ \lambda=(1,1,2,3,4,4,5,5,6,6,8,9)\]
 can be expressed as follows
\begin{equation} \label{gm}
\lambda=\setcounter{MaxMatrixCols}{9}\begin{bmatrix}
\ &\ &\ &\ &5&\ &\ &\ &\ \\[3pt]
 \ &2&\ &4&\ &6 &\ &\ &\ \\[3pt]
1&\ &\ &4&\ &6&\ &\ &9\\[3pt]
1&\ &3&\ &5&\ &\ &8&\
\end{bmatrix}\begin{matrix}
4\\[3pt]3\\[3pt]2\\[3pt]1
\end{matrix} \; ,
\end{equation}
where the marks are listed outside the brackets, that is,
the parts at the bottom are marked with $1$, and the parts immediately next to the
bottom line are marked by $2$, and so on. The Gordon marking of a partition can be
 considered as a way to represent a partition. For this reason, the diagram (\ref{gm})
 is called the Gordon marking representation of a partition.

We shall introduce the Gordon marking of an overpartition. In fact,
the three bijections in the proof of Theorem \ref{main}
are constructed  based on Gordon markings of overpartitions.
 The Gordon marking of
an overpartition can be defined as follows. It is clear that this notion is an extension of
the Gordon marking of an ordinary  partition.

\begin{defi}
The Gordon marking of an overpartition $\lambda$ is an assignment of
positive integers (marks) to parts of $\lambda$. We assign the marks to parts
in the following order
\begin{equation}\label{order}
 \overline{1} < 1 < \overline{2} < 2< \cdots
 \end{equation}
such that the marks are as small as possible subject to the following conditions.
If $\overline{j+1}$ is not a part of $\lambda$,
 then all the parts $j$, $\overline{j}$, and $j+1$ are assigned  different integers.
  If  $\lambda$ contains an overlined part $\overline{j+1}$,
   then the smallest mark
assigned to a part $j$ or $\overline{j}$ can be used as the mark of  $j+1$ or $\overline{j+1}$.
\end{defi}

For example, given an overpartition
\[ \lambda= (16,13,12,12,11,\overline{10},\overline{8},8,8,7,\overline{6},6,5,5,4,2,2,\overline{1}).\]
 The
Gordon marking of $\lambda$ is
\[ (\overline{1}_1,2_2,2_3,4_1,5_2,5_3,
\overline{6}_1,6_2,7_3,\overline{8}_1,8_2,8_3,\overline{10}_1,11_2,12_1,12_3,13_2,16_1),
\]
where the subscripts are the marks. The Gordon marking of $\lambda$ can
also be illustrated as
\begin{equation}\label{lambda}
\lambda=\setcounter{MaxMatrixCols}{16}\begin{bmatrix}
 \ &2&\ &\ &5&\ &7&8&\ &\ &\ &12&\ &\ &\ &\\[3pt]
\ &2&\ &\ &5&6&\ &8&\ &\ &11&\ &13&\ &\ &\\[3pt]
\overline{1}&\ &\ &4&\ &\overline{6}&\ &\overline{8}&\ &\overline{10}&\ &12&\ &\ &\ &16
\end{bmatrix}\begin{matrix}
3\\[3pt]2\\[3pt]1
\end{matrix}\;,
\end{equation}
where
the parts in the third  row are marked by $1$, the parts in the second row
are marked by $2$, and the parts in the first row are marked by $3$.

It is not hard to see that the Gordon marking of any overpartition is unique.
To compute the generating function of $D_{k,i}(m,n)$,  let  $T_{k,i}(m,n)$
denote the set of overpartitions enumerated by $D_{k,i}(m,n)$.
 We further classify   $T_{k,i}(m,n)$ by considering whether
 the smallest part of an overpartition is overlined element. Keep in mind that
 the parts of an overpartition are ordered by \eqref{order}.
 Let $U_{k,i}(m,n)$
 denote the set of overpartitions
in $T_{k,i}(m,n)$ for which  the smallest
part is overlined,  and let $I_{k,i}(m,n)$ denote the set of overpartitions in $T_{k,i}(m,n)$
with non-overlined smallest part.
Thus we have
 \begin{equation} \label{uki}
  T_{k,i}(m,n)= U_{k,i}(m,n)  \cup I_{k,i}(m,n).
  \end{equation}
Let $F_{k,i}(m,n)=|U_{k,i}(m,n)|$ and $G_{k,i}(m,n)=|I_{k,i}(m,n)|$. Then we have
\begin{equation}\label{DFG}D_{k,i}(m,n)=F_{k,i}(m,n)+G_{k,i}(m,n).\end{equation}

 Below is a relation between $F_{k,i}(m,n)$ and $G_{k,i}(m,n)$.

 \begin{lem}\label{FG1}For $2\leq i\leq k$, we have
 \begin{equation}\label{FG} F_{k,i-1}(m,n)=G_{k,i}(m,n).\end{equation}
 For $i=1$, we  have \begin{equation}\label{FG2}G_{k,1}(m,n)=F_{k,k}(m,n-m).\end{equation}
 \end{lem}

\pf For $i\geq 2$, there is a simple bijection between
$U_{k,i-1}(m,n)$ and $I_{k,i}(m,n)$.
For an overpartition $\lambda\in U_{k,i-1}(m,n)$, we change the smallest  part $\overline{j}$  of $\lambda$ to a
non-overlined part $j$. Then we get an overpartition in $I_{k,i}(m,n)$.
Conversely, we can change one of  the smallest part $j$ of
an overpartition $\beta \in I_{k,i}(m,n) $  to an overlined part $\overline{j}$
 to get an overpartition in $U_{k,i-1}(m,n)$.
Clearly, this map is a bijection. Hence  \eqref{FG} holds for $i\geq 2$.

For $i=1$, we shall show a bijection between  $I_{k,1}(m,n)$ and  $U_{k,k}(m,n-m)$.
 Substracting one from each part of overpartition $\lambda$ in $I_{k,1}(m,n)$ and changing
  one of the smallest parts to an overlined part, we obtain  an overpartition in $U_{k,k}(m,n-m)$.
Conversely, for an overpartition $\mu$ in $U_{k,k}(m,n-m)$, we can switch
the smallest part to a non-overlined part, and increase each part of $\mu$ by one (regardless
of the overlines), so that we can get an overpartition in $I_{k,1}(m,n)$.
So we arrive at \eqref{FG2}. This completes the proof.  \qed

By the above lemma,  the generating function of $G_{k,i}(m,n)$ can be obtained
 from the generating function of  $F_{k,i}(m,n)$. Moreover, from  \eqref{DFG} it follows that
  the generating function of $D_{k,i}(m,n)$ can be deduced from $F_{k,i}(m,n)$.
The following theorem gives the generating function of $F_{k,i}(m,n)$.

\begin{thm}\label{F}
For $k\geq i\geq 1$,
\begin{align}
&\sum_{n=0}^{\infty}F_{k,i}(m,n)x^mq^n\nonumber\\&\qquad
=\sum_{N_1\geq N_2\geq\cdots\geq
N_{k-1}\geq0}\frac{q^{\frac{(N_1+1)N_1}{2}+N_2^2+\cdots+N_{k-1}^2+N_{i+1}+\cdots+N_{k-1}}
(-q)_{N_1-1}x^{N_1+\cdots+N_{k-1}}}{(q)_{N_1-N_2}\cdots(q)_{N_{k-2}-N_{k-1}}(q)_{N_{k-1}}}.\end{align}
 \end{thm}

To derive the generating function of
  $F_{k,i}(m,n)$, we shall further classify the set
   $U_{k,i}(m,n)$.
Let $\lambda^{(r)}$ denote the
partition that consists of all $r$-marked parts of $\lambda$.
Let $N_r$ be the number of $r$-marked parts (i.e. the number of
parts in $\lambda^{(r)}$), and let $n_r=N_r-N_{r-1}$ for any positive
integer $r$. Notice that for any overpartition $\lambda$ enumerated by $D_{k,i}(m,n)$,
 the parts $j$, $\overline{j}$ and  $j+1$ occur at most $k-1$ times in
$\lambda$. It follows that the marks of $\lambda$ do not exceed $k-1$.
So we are led to consider the  parameters $N_1,\ldots, N_{k-1}$ and $n_1,\ldots,n_{k-1}$
as the summation indices when we compute the generating function of $F_{k,i}(m,n)$.
 It also can be seen that $N_1\geq N_2\geq \cdots\geq N_{k-1}\geq 0$ and  $n_1, n_2,
\ldots, n_{k-1}\geq0$.
The detailed proof of Theorem \ref{F} will be given
in the next four sections.

\section{The first bijection for the proof of Theorem \ref{main}}

In this section, we  classify the set $U_{k,i}(m,n)$ according to the parameters $N_1,\ldots,N_{k-1}$, and we give the first bijection for the proof of Theorem \ref{main}. Let $\sum_{i=1}^{k-1}N_i=m$, and
let $U_{N_1,N_2,\ldots,N_{k-1};i}(n)$ denote the set of overpartitions   in $U_{k,i}(m,n)$
that have  $N_r$ $r$-marked parts for $1\leq r\leq k-1$. Let $P_{N_1,N_2,\ldots,N_{k-1};i}(n)$ denote the set of overpartitions in  $U_{N_1,N_2,\ldots,N_{k-1};i}(n)$ for which all the 1-marked parts are overlined. Set
\begin{align}&U_{N_1,N_2,\ldots,N_{k-1};i}=\bigcup_{n\geq 0}U_{N_1,N_2,\ldots,N_{k-1};i}(n),\\
&P_{N_1,N_2,\ldots,N_{k-1};i}=\bigcup_{n\geq 0}P_{N_1,N_2,\ldots,N_{k-1};i}(n).
\end{align}

More precisely, we shall give a bijection
for the following relation.

\begin{thm} \label{FP}For $k\geq i \geq 1$, we have
\begin{equation}
  \sum_{\lambda\in
U_{N_1,N_2,\cdots,N_{k-1};i}}x^{l(\lambda)}q^{|\lambda|}=(-q)_{N_1-1}
\sum_{\alpha\in
P_{N_1,N_2,\cdots,N_{k-1};i}}x^{l(\alpha)}q^{|\alpha|},
\end{equation}
where  $l(\lambda)$ denotes the number of parts of $\lambda$.
\end{thm}

 Before we present the bijection for the above relation, we introduce a reduction  operation
 based on the Gordon markings, which transforms an overpartition in $U_{N_1,N_2,\ldots,N_{k-1};i}(n)$
 containing at least one non-overlined part with mark 1 to
an overpartition in $U_{N_1,N_2,\ldots,N_{k-1};i}(n-1)$. This reduction operation
preserves the number of $r$-marked parts for $r=1,2,\ldots, k-1$. Since we shall give another
reduction operation in the next section, we call the reduction operation described below
the first reduction operation.

\noindent
{\bf The First Reduction Operation.}
Let $\lambda=(\lambda_1,\ldots,\lambda_m)$ be an
overpartition of $n$ containing at least one non-overlined part with mark $1$.
Assume that $\lambda_j$ is the rightmost
 non-overlined part with mark 1. To be more precisely, for a part $\lambda_j$,
we write $\lambda_j=\overline{a_j}$ to indicate that $\lambda_j$ is an overline part
and write $\lambda_j=a_j$ to indicate that $\lambda_j$ is a non-overline part. Moreover,
we say that $a_j$ is the underlying part of $\lambda_j$. We consider two cases.

\noindent
Case 1.  There is a non-overlined part $a_j+1$ of $\lambda$  but there is no overlined $1$-marked part $\overline{a_j+1}$.
First, we change the part $\lambda_j$ to a $1$-mark  part $\overline{a_j}$. Then we choose
the part  $a_j+1$ with the smallest mark, say $r$, and replace this
 $r$-marked
 part $a_j+1$ with a $r$-marked part $a_j$.
Since in $\lambda$ $r$ is the smallest mark of the parts $a_j+1$ and the
$1$-marked $a_j$ is non-overlined, by the definition of the
 Gordon marking of an  overpartition, we deduce
  that  either $r$ is still the smallest mark of the parts with underlying part $a_j-1$ or there are no parts
 with underlying part $a_j-1$. In either case, we may place the new $r$-marked part $a_j$ in a position with mark $r$.

If there is a $1$-marked overlined part to the right of the $\overline{a_j}$, we switch it to a non-overlined part and we can see that the rightmost $1$-marked nonoverlined part of the
resulting overpartition is right to $\lambda_j$. If there are no $1$-marked parts larger than $a_j$,
we shall do nothing and in this case we can notice that the number of $1$-marked overlined parts in the resulting
overpartition is one more than it in $\lambda$.
  In either case, we denote the resulting
 overpartition by $\mu$. Clearly, $\mu$ is an overpartition of $n-1$. Moreover,
 it can be seen that $\mu$
 contains the same number of $r$-marked parts as $\lambda$, for $1\leq r\leq k-1$.

\noindent
 Case 2. Either  an overlined part $\overline{a_j+1}$ is a $1$-marked  part of $\lambda$
  or  there are no  parts with underlying part $a_j+1$. In either case,
    we may change the part $\lambda_j$ to a $1$-marked  overlined part $\overline{a_j-1}$.

If  there are $1$-marked parts larger than $a_j$,
then they are all overlined parts because of the choice of $\lambda_j$.
In this case we switch the overlined $1$-marked part next to $\lambda_j$
to a non-overlined part. Let $\mu$ denote the resulting overpartition.
 It is easily seen that
 in this case  the rightmost non-overlined part in $\mu$ is right to the part
 $\lambda_j$ and $\mu$ has  the same number of $1$-marked overlined parts and the same number of $1$-marked  nonoverlined parts as $\lambda$.

 It remains to consider the
 case when there are no $1$-marked parts larger than $a_j$.
 In this case, no operation is needed and we  set $\mu$ to be the   overpartition obtained in the
 previous step.
 It is clear that
 $\mu$ has one more $1$-marked overlined parts and  one less $1$-marked non-overlined parts than $\lambda$.

  In either case, one can deduce that $\mu$ is an overpartition
of $n-1$ with the same number of $r$-marked parts as $\lambda$, for $1\leq r \leq k-1$.

For example, let  $\lambda$ be an overpartition  in $U_{7,6,5;1}(135)$ as given below
\begin{equation}
\setcounter{MaxMatrixCols}{15}\begin{bmatrix}\nonumber
 \ &2&\ &\ &5&\ &7&8&\ &\ &\ &12&\ &\ &\ \\[3pt]
\ &2&\ &\ &5&6&\ &8&\ &\ &11&\ &\mathbf{13}&\ &\ \\[3pt]
\overline{1}&\ &\ &4&\ &\overline{6}&\ &\overline{8}&\ &\overline{10}&\ &\mathbf{12}&\ &\ &\mathbf{\overline{15}}
\end{bmatrix}\begin{matrix}
3\\[3pt]2\\[3pt]1
\end{matrix}\;.
\end{equation}
The part ${\bf 12}$ with mark $1$ is the $\lambda_j$ as in the
 description of the reduction operation,
 since it is the rightmost non-overlined part with mark $1$.
 Notice that $13$ is not a $1$-marked part of $\lambda$,
         but ${\bf 13}$ is a $2$-marked part. By the operation in Case 1, we change the
         $1$-marked part ${\bf 12}$
 to a part ${\bf\overline{12}}$, then we change the
 $2$-marked
 part ${\bf 13}$ to   ${\bf 12}$ and place it in a position with mark $2$.
 Then we switch ${\bf \overline{15}}$  to ${\bf 15}$.
After the reduction operation by choose $\lambda_j$ to be $1$-marked $12$, we get an overpartition $\mu$ in
$U_{7,6,5;1}(134)$
\begin{equation}
\setcounter{MaxMatrixCols}{15}\begin{bmatrix}\nonumber
 \ &2&\ &\ &5&\ &7&8&\ &\ &\ &12&\ &\ &\ \\[3pt]
\ &2&\ &\ &5&6&\ &8&\ &\ &11&\mathbf{12} &\ &\ &\ \\[3pt]
\overline{1}&\ &\ &4&\ &\overline{6}&\ &\overline{8}&\ &\overline{10}&\ &\mathbf{\overline{12}}&\ &\ &\mathbf{15}
\end{bmatrix}\begin{matrix}
3\\[3pt]2\\[3pt]1
\end{matrix}\;.
\end{equation}

 Let us apply the reduction operation to above overpartition $\mu$.
 The part ${\bf 15}$ is the rightmost non-overlined part with mark $1$ in $\mu$ and
  there are no parts greater than $15$. So we need to apply the operation in  Case 2.
  By changing  ${\bf 15}$ to ${\bf \overline{14}}$, we obtain an  overpartition in $U_{7,6,5;1}(133)$
\begin{equation}
\setcounter{MaxMatrixCols}{14}\begin{bmatrix}\nonumber
 \ &2&\ &\ &5&\ &7&8&\ &\ &\ &12&\ &\  \\[3pt]
\ &2&\ &\ &5&6&\ &8&\ &\ &11&12 &\ &\  \\[3pt]
\overline{1}&\ &\ &4&\ &\overline{6}&\ &\overline{8}&\ &\overline{10}&\ &\overline{12}&\ &\mathbf{\overline{14}}
\end{bmatrix}\begin{matrix}
3\\[3pt]2\\[3pt]1
\end{matrix}\;.
\end{equation}

Indeed, the above reduction operation is reversible. This implies that
there is a bijection for the relation in  Theorem 4.1.
We shall give the  dilation operation as the inverse of the reduction operation,
and we shall call it the first dilation operation. In fact, there are two types of
dilation operations depending on the choice of the position where
the operation will take place.

\noindent
{\bf The First Dilation Operation.}
Let $\lambda=(\lambda_1,\ldots,\lambda_m)$ be an
overpartition in $U_{N_1,N_2,\ldots,N_{k-1};i}(n)$. For a part $\lambda_j$,
we use $a_j$ to denote the underlying part of $\lambda_j$.

We proceed to determine the   part $\lambda_j$
 which tells where  the dilation operation will take place. There are two types of the dilation operation.
If there are no $1$-marked parts next to the rightmost overlined part $\lambda_j$,
then we may choose $\lambda_j$ and we shall say that  the operation is of type $A$.
If there is at least one overlined part such that
the next $1$-marked part is  non-overlined, then we choose the rightmost one to be $\lambda_j$.
For this choice,
 we   say that the dilation operation is of type $B$.
 It should be mentioned that it is possible that we can apply two types of
 operations to an overpartition. For
 each overpartition in $U_{N_1,N_2,\ldots,N_{k-1};i}(n)$,
 we can apply at least one of the two types of the dilation operation.
 As will be seen, in the proof of Theorem 4.1 we need to consider how to apply
 the two types of the dilation operation.

\noindent
Case 1:  There are two parts of the same mark with underlying parts $a_j$ and $a_j-1$, we denote this same mark by $r$.
 It should be noticed that there are no $1$-marked
 parts  with underlying part $a_j+1$ because of the choice of $\lambda_j$.
  We change $\lambda_j$ to a non-overlined part $a_j$ and replace the $r$-marked
  part   $a_j$ by an  $r$-marked  part  $a_j+1$.

If there are $1$-marked parts with underlying parts greater than $a_j$, we consider
the leftmost one, which must be non-overlined, and we change this
 non-overlined part to an  $1$-marked  overlined part. Denote the resulting overpartition by $\mu$. Clearly, the rightmost $1$-marked overlined part to the left of a non-overlined part in $\mu$  must be to the left of
 $\lambda_j$ in $\lambda$. Moreover,
   $\mu$ has  the same number of $1$-marked overlined parts and the same number of $1$-marked
 non-overlined parts as $\lambda$.

 We now turn to the case when there are no $1$-marked parts with underlying parts greater
than $a_j$. In this case no operation is required and we
denote the  overpartition obtained so far by $\mu$.
Notice that $\mu$ has one less $1$-marked overlined parts and  one more $1$-marked non-overlined parts than $\lambda$.

In either case, one can deduce that $\mu$ is an overpartition
in $U_{N_1,N_2,\ldots,N_{k-1};i}(n+1)$ with the same number of $r$-marked parts as $\lambda$,
for $1\leq r \leq k-1$.

\noindent
Case 2:
There are no two parts  with underlying parts $a_j$ and $a_j-1$ that have the same mark.
We see that there is no $1$-marked part with underlying part $a_j+1$ because of the choice of
$\lambda_j$. We
 change $\lambda_j$ to  a non-overlined part $a_j$ with mark $1$.
We   denote by $r$
 the largest mark  of the parts equal to $a_j$, and replace the
$r$-marked non-overlined part $a_j$ with an $r$-marked non-overlined part $a_j+1$.
Since $r$ is the largest mark of the parts equal to $a_j$ and $a_j+1$ is not a $1$-marked part  of $\lambda$, we see that $a_j+1$  cannot be a part with a mark not exceeding $r$.
So we may place the new part equal to  $a_j+1$ in a position of mark $r$.

If there is a  $1$-marked non-overlined  part next to $\lambda_j$,
we switch this non-overlined part
to an overlined part.
Let $\mu$ denote the resulting overpartition.
It is easily seen that in this case
$\mu$ has  the same number of $1$-marked overlined parts and the
same number of $1$-marked  non-overlined parts as $\lambda$.

We still need to consider the
 case when  there are no parts next to $\lambda_j$,
 In this case, we just denote the resulting overpartition by $\mu$.
 Clearly,  $\mu$ has one more $1$-marked non-overlined parts and  one less $1$-marked overlined parts than $\lambda$.

In either case, we see that $\mu$ is an overpartition
 in $U_{N_1,N_2,\ldots,N_{k-1};i}(n+1)$ with the same number of $r$-marked parts as $\lambda$,
for $1\leq r \leq k-1$.

It is easily checked that the first reduction operation  is the inverse
of the first dilation operation. More precisely, we have the following
property.

\begin{thm}
The dilation operation of Type A is the inverse of the reduction operation which increases the number of overlined parts in  $\lambda$, whereas the dilation operation of Type B is the inverse of the reduction operation which preserves the number of overlined parts in $\lambda$.
\end{thm}

We are now ready to present the proof of Theorem \ref{FP}.

\noindent{\it Proof of Theorem \ref{FP}.}
Based on the reduction operation, we shall establish  a bijection $\varphi$ between
$U_{N_1,N_2,\ldots,N_{k-1};i}$ and $P_{N_1,N_2,\ldots,N_{k-1};i}\times D_{N_1}$, where $D_{N_1}$ denotes the set of ordinary partitions with distinct parts such that
 each part is less than $N_1$.
Let $\lambda$ be an overpartition in $ U_{N_1,N_2,\ldots,N_{k-1};i}$. We shall give
a procedure to construct $\varphi(\lambda)$, which is
a pair $(\alpha, \beta)$, where $\alpha$ is an overpartition in $P_{N_1,N_2,\ldots,N_{k-1};i}$ and $\beta$ is a partition in $D_{N_1}$.

\noindent Step 1.  Set $\alpha=\lambda$, $\beta=\phi$ and $t=1$. If there are no non-overlined $1$-marked parts in $\alpha$, go to Step 3; Otherwise, go to Step 2.

\noindent Step 2. If the largest $1$-marked part of $\alpha$ is overlined,
then apply the first reduction operation on  $\alpha$.  If there are still
 non-overlined $1$-marked  parts in $\alpha$,
then set  $t$ to $t+1$ and repeat this step; Otherwise, go to Step 3.

If  the largest $1$-marked part of $\alpha$ is non-overlined, then add $t$ to $\beta$ as a
new part and apply the first reduction operation on $\alpha$.
Reset $t$ to $1$ and repeat this step if there are still non-overlined $1$-marked  parts in $\alpha$;
Otherwise, go to Step 3.

\noindent Step 3.  Set $\varphi(\lambda)=(\alpha, \beta)$.

Evidently, $\alpha$ is an overpartition  in $P_{N_1,N_2,\ldots,N_{k-1};i}$ and $|\lambda|=|\alpha|+|\beta|$. It remains to prove that the parts of
$\beta$ are less than $N_1$.
Let \[\lambda^{(1)}_1 < \lambda^{(1)}_2<\cdots<\lambda^{(1)}_{N_1}\] denote the $1$-marked parts of $\lambda$. Moreover,
suppose that there are $s$ non-overlined $1$-marked parts of $\lambda$, which are
denoted by
\[ \lambda^{(1)}_{i_{1}}  <\lambda^{(1)}_{i_{2}}  <\cdots<\lambda^{(1)}_{i_{s}} .\]
Examining Step 2 of the above procedure, we see that
after applying the operation in Step
      2 to the  rightmost non-overlined part such that it is the largest $1$-marked part of $\alpha$, the number of non-overlined part decreases by one. So we find that
 for each non-overlined $1$-marked part $\lambda_{i_{t}}^{(1)}$, we can
    iterate Step 2 $N_1-i_t+1$ times in order to decrease the number of non-overlined parts
    by one and add $N_1-i_t+1$ to $\beta$ as a new part.
    Hence  we deduce that
$\beta=(N_1-i_1+1, N_1-i_2+1, \ldots, N_1-i_s+1)$.
Recall that the smallest $1$-marked part of an overpartition in $U_{N_1,N_2,\ldots,N_{k-1};i}$ is always overlined. It follows
that $N_1-i_t+1<N_1$, for $1\leq t\leq s$. So $\beta$ is a partition in $D_{N_1}$.

Next we give a brief description of the inverse of  $\varphi$. The detailed proof is
omitted because it is a straightforward verification.

Let $\alpha$ be an overpartition in $P_{N_1,N_2,\ldots N_k;i}$ and
$\beta=(\beta_1,\beta_2,\cdots,\beta_s)$
 be a  partition with distinct parts and  $\beta_1\leq N_1-1$.
 We shall give
a procedure to construct $\varphi^{-1}(\alpha,\beta)$, which is
an overpatition $\lambda$ in $U_{N_1,N_2,\ldots N_k;i}$.

\noindent Step 1. Set $\lambda=\alpha$. Let $s$ be the number of parts in $\beta$.

\noindent Step 2.
For $t$ from $1$ to $s$,
apply the  dilation operation  of type $A$  to $\lambda$. Then
the   dilation operation of type $B$ will be applied $\beta_t-1$ times  to   $\lambda$.
Now we get an overpartion $\lambda$ in $U_{N_1,N_2,\ldots N_k;i}$.
It can be checked  that
$\lambda^{(1)}_{N_1-\beta_1},\ldots,\lambda^{(1)}_{N_1-\beta_s}$  are the non-overlined $1$-marked
parts of $\lambda$.

To prove that $\varphi^{-1}(\varphi(\lambda))=\lambda$,
 we need the fact that the first reduction operation
and the first dilation operation are inverses of each other.
This completes the proof. \qed

To demonstrate the above bijection we give an example. Let $\lambda$ be the overpartition as given in \eqref{lambda}, that is,

\[
\lambda=\setcounter{MaxMatrixCols}{16}\begin{bmatrix}
 \ &2&\ &\ &5&\ &7&8&\ &\ &\ &12&\ &\ &\ &\\[3pt]
\ &2&\ &\ &5&6&\ &8&\ &\ &11&\ &13&\ &\ &\\[3pt]
\overline{1}&\ &\ &4&\ &\overline{6}&\ &\overline{8}&\ &\overline{10}&\ &12&\ &\ &\ &16
\end{bmatrix}\begin{matrix}
3\\[3pt]2\\[3pt]1
\end{matrix}\;.
\]

First we set $\alpha=\lambda$, $\beta=\phi$ and $t=1$.
Notice that the greatest  $1$-marked part of $\alpha$ is non-overlined. So we
 let $\beta=(1)$ and set $t=1$.
Applying the first reduction operation,  we have
\begin{equation}
\setcounter{MaxMatrixCols}{15}\alpha=\begin{bmatrix}\nonumber
 \ &2&\ &\ &5&\ &7&8&\ &\ &\ &12&\ &\ &\ \\[3pt]
\ &2&\ &\ &5&6&\ &8&\ &\ &11&\ &13&\ &\ \\[3pt]
\overline{1}&\ &\ &4&\ &\overline{6}&\ &\overline{8}&\ &\overline{10}&\ &12&\ &\ &\overline{15}
\end{bmatrix}\begin{matrix}
3\\[3pt]2\\[3pt]1
\end{matrix}\;.
\end{equation}
Since the greatest non-overlined $1$-marked part is  $12$  which  is not the greatest $1$-marked part, we apply the first reduction operation on $\alpha$ and  let $t=2$.
Then we get
\begin{equation}
\setcounter{MaxMatrixCols}{15}\alpha=\begin{bmatrix}\nonumber
 \ &2&\ &\ &5&\ &7&8&\ &\ &\ &12&\ &\ &\ \\[3pt]
\ &2&\ &\ &5&6&\ &8&\ &\ &11&12 &\ &\ &\ \\[3pt]
\overline{1}&\ &\ &4&\ &\overline{6}&\ &\overline{8}&\ &\overline{10}&\ &\overline{12}&\ &\ &15
\end{bmatrix}\begin{matrix}
3\\[3pt]2\\[3pt]1
\end{matrix}\;.
\end{equation}
Now the rightmost non-overlined $1$-marked part is $15$ and it is the greatest $1$-marked part.
So we
apply the reduction operation and  let $\beta=(2,1)$.  Now we should  reset $t=1$. Then we get
\begin{equation}
\setcounter{MaxMatrixCols}{14}\alpha=\begin{bmatrix}\nonumber
 \ &2&\ &\ &5&\ &7&8&\ &\ &\ &12&\ &\  \\[3pt]
\ &2&\ &\ &5&6&\ &8&\ &\ &11&12 &\ &\  \\[3pt]
\overline{1}&\ &\ &4&\ &\overline{6}&\ &\overline{8}&\ &\overline{10}&\ &\overline{12}&\ &\overline{14}
\end{bmatrix}\begin{matrix}
3\\[3pt]2\\[3pt]1
\end{matrix}\;.
\end{equation}
In order to get an overpartition with no non-overlined $1$-marked parts, we still need to apply the reduction operation $6$ times. The details are omitted. Finally, we obtain
\begin{equation}\label{alpha}
\setcounter{MaxMatrixCols}{13}\alpha=\begin{bmatrix}
 \ &2&\ &\ &5&6 &\ &8&\ &\ &\ &12&\\[3pt]
\ &2&\ &4 &\ &6&\ &8&\ &10 &\ &12 &\\[3pt]
\overline{1}&\ &\ &\overline{4}&\ &\overline{6}&\overline{7}&\ &\ &\overline{10}&\overline{11} &\ &\overline{13}
\end{bmatrix}\begin{matrix}
3\\[3pt]2\\[3pt]1
\end{matrix}\;,
\end{equation}
and  $\beta=(6,3,1)$. Thus we have constructed a pair $(\alpha, \beta)$,
 where $\alpha$ is an overpartition such that all $1$-marked parts
overlined, $\beta$ is partition in $D_{7}$. Moreover, we have
$|\lambda|=|\alpha|+|\beta|$.

\section{The second bijection for the proof of Theorem \ref{main}}

In this section, we  introduce a class of overpartitions in
$P_{N_1,N_2,\ldots,N_{k-1};i}$, which will be denoted by $Q_{N_1,N_2,\ldots,N_{k-1};i}$.
We aim to relate the generating function of $P_{N_1,N_2,\ldots,N_{k-1};i}$ to that of
$Q_{N_1,N_2,\ldots,N_{k-1};i}$. To define the set $Q_{N_1,N_2,\ldots,N_{k-1};i}$, we observe that for any $\lambda \in P_{N_1,N_2,\ldots,N_{k-1};i}(n)$ and for any $1\leq t\leq n$, we have
 \begin{equation}\label{ft}
 f_t(\lambda)+f_{\overline{t}}(\lambda)+f_{t+1}(\lambda)\leq  k-1,
  \end{equation}
  where $f_t(\lambda)$ denotes the number of occurrences of $t$ in $\lambda$.
 We define the set $Q_{N_1,N_2,\ldots,N_{k-1};i}$ as the set of overparitions $\lambda$ in $P_{N_1,N_2,\ldots,N_{k-1};i}$ for which the equality holds in (\ref{ft}), namely,
 \begin{equation}\label{fte}
 f_t(\lambda)+f_{\overline{t}}(\lambda)+f_{t+1}(\lambda) =  k-1
  \end{equation}
  for any positive integer $t$ which is smaller than the greatest $(k-1)$-marked part.
 It should be mentioned that  Bressoud \cite{Bre79, Bre80} obtained a generalization of the Rogers-Ramanujan identities by considering ordinary partitions $\lambda$
  that satisfy the equality in (\ref{fte}), namely,
  \begin{equation}
  f_t(\lambda)+f_{t+1}(\lambda) = k-1.
  \end{equation}

Set
\[Q_{N_1,N_2,\ldots,N_{k-1};i}=\bigcup_{n\geq 0}Q_{N_1,N_2,\ldots,N_{k-1};i}(n).\]

The following theorem  establishes
a relation between the generating function of $P_{N_1,N_2,\ldots,N_{k-1};i}$ and the generating function of $Q_{N_1,N_2,\ldots,N_{k-1};i}$.

\begin{thm}\label{PQ} For $N_1\geq N_2 \geq \cdots\geq N_{k-1}\geq 0$, we have
 \begin{equation}
  \sum_{\alpha\in P_{N_1,N_2,\ldots,N_{k-1};i}}x^{l(\alpha)}q^{|\alpha|}
=\frac{1}{(q)_{N_{k-1}}}\sum_{\gamma\in Q_{N_1,N_2,\ldots,N_{k-1};i}}x^{l(\gamma)}q^{|\gamma|}.
 \end{equation}
\end{thm}

To prove the above theorem, we shall give a bijection  based on
a reduction operation and a dilation operation which are called
 the second reduction  and the second dilation.
The second reduction transforms an overpartition  $\alpha$ in
$P_{N_1,N_2,\ldots,N_{k-1};i}(n)\setminus Q_{N_1,N_2,\ldots,N_{k-1};i}(n)$ to an
overpartition in $P_{N_1,N_2,\ldots,N_{k-1};i}(n-1)$. More precisely,
this operation requires the choice of a $(k-1)$-marked part $\alpha_j$ whose
underlying part is $t$ satisfying one of the following two conditions

\begin{itemize}
\item[1.] There are no parts with underlying part  $t-1$;
\item[2.] There is a part with underlying part  $t-1$ and
\begin{equation} \label{ft3}
 f_{t-2}(\alpha)+f_{\overline{t-2}}(\alpha)+f_{t-1}(\alpha)<k-1.
  \end{equation}
\end{itemize}

By the definitions of $P_{N_1,N_2,\ldots,N_{k-1};i}(n)$ and $Q_{N_1,N_2,\ldots,N_{k-1};i}(n-1)$, it is not difficult to see that for any $\alpha$ in $P_{N_1,N_2,\ldots,N_{k-1};i}(n)\setminus Q_{N_1,N_2,\ldots,N_{k-1};i}(n)$, there
 exists a $(k-1)$-marked part  $\alpha_j$ satisfying one of the above conditions.

\noindent{\bf The Second Reduction Operation.}
Let $\alpha=(\alpha_1,\ldots,\alpha_m)$ be an
overpartition in $P_{N_1,N_2,\ldots,N_{k-1};i}(n)\setminus Q_{N_1,N_2,\ldots,N_{k-1};i}(n)$.
Let  $\alpha_j$ be a  $(k-1)$-marked part
with underlying part $t$ satisfying  one of the above conditions.

If $\alpha_j$ satisfies Condition 1, that is, there are no parts with underlying part $t-1$,
then there is an overlined part $\overline{t}$ since $t$ is the underlying part of $\alpha_j$.
 We replace  $\overline{t}$ with a  $1$-marked overlined part $\overline{t-1}$.

If $\alpha_j$ satisfies Condition 2, write (\ref{ft3}) as
\begin{equation}\label{ft4}
 \sum_{l=1}^{k-1}(f_{t-2}(\alpha^{(l)})+f_{\overline{t-2}}(\alpha^{(l)})+f_{t-1}(\alpha^{(l)}))<k-1,
  \end{equation}
  where $\alpha^{(l)}$ is the overpartition consisting of the $l$-marked parts of $\alpha$.
 So we can find the smallest mark $r\geq 2$ such that $t$ is a part of mark $r$ and
  \begin{equation}\label{ftc}
 \sum_{l=1}^{r}(f_{t-2}(\alpha^{(l)})+f_{\overline{t-2}}(\alpha^{(l)})+f_{t-1}(\alpha^{(l)}))<r.
  \end{equation}
  Replace  the $r$-marked part $t$ with an $r$-marked part  $t-1$.

It can be seen that in either case we obtain  the Gordon
 marking representation of an overpartition  in $P_{N_1,N_2,\ldots,N_{k-1};i}(n-1)$.

For example, let  $\alpha$ be an overpartition  in $P_{7,6,5;1}(126)$ as given below
\begin{equation}
\setcounter{MaxMatrixCols}{13}\alpha=\begin{bmatrix}\nonumber
 \ &2&\ &4&\ &6 &\ &8&\ &\ &\ &12&\\[3pt]
\ &2&\ &4 &\ &6&\ &8&\ &10 &\ &12 &\\[3pt]
\overline{1}&\ &\ &{\bf \overline{4}}&\ &\overline{6}&\overline{7}&\ &\ &\overline{10}&\overline{11} &\ &\overline{13}
\end{bmatrix}\begin{matrix}
3\\[3pt]2\\[3pt]1
\end{matrix}\; .
\end{equation}
Choosing $\alpha_j$ to be the $3$-marked part $4$,  we see that the $1$-marked part $\overline{4}$ satisfies
 Condition  1. Then we can replace $\overline{4}$ with a $1$-marked $\overline{3}$ to transform $\alpha$ to
 an overpartition in $P_{7,6,5;1}(125)$:
\begin{equation}
\setcounter{MaxMatrixCols}{13}\begin{bmatrix}\nonumber
 \ &2&\ &4&\ &6 &\ &8&\ &\ &\ &12&\\[3pt]
\ &2&\ &{\bf 4} &\ &6&\ &8&\ &10 &\ &12 &\\[3pt]
\overline{1}&\ &\overline{3} &\ &\ &\overline{6}&\overline{7}&\ &\ &\overline{10}&\overline{11} &\ &\overline{13}
\end{bmatrix}\begin{matrix}
3\\[3pt]2\\[3pt]1
\end{matrix}\;.
\end{equation}

For the above overpartition, choosing the same $\alpha_j$ as before,  we
 see that the $3$-marked part $4$ satisfies Condition 2.
We  further apply the reduction in this case.
 Clearly,   $2$ is the smallest mark satisfying   Condition \eqref{ftc}.
 So we can replace the $2$-marked part $4$ with a $2$-marked part $3$ to form an overpartition in $P_{7,6,5;1}(124)$:
\begin{equation}
\setcounter{MaxMatrixCols}{13}\begin{bmatrix}\nonumber
 \ &2&\ &4&\ &6 &\ &8&\ &\ &\ &12&\\[3pt]
\ &2&3 &\ &\ &6&\ &8&\ &10 &\ &12 &\\[3pt]
\overline{1}&\ &\overline{3} &\ &\ &\overline{6}&\overline{7}&\ &\ &\overline{10}&\overline{11} &\ &\overline{13}
\end{bmatrix}\begin{matrix}
3\\[3pt]2\\[3pt]1
\end{matrix}\;.
\end{equation}

The second dilation transforms an overpartition  $\alpha$ in
$P_{N_1,N_2,\ldots,N_{k-1};i}(n)$ to an overpartition in $P_{N_1,N_2,\ldots,N_{k-1};i}(n-1)\setminus Q_{N_1,N_2,\ldots,N_{k-1};i}(n-1)$.
To be more specific, the operation starts with a choice of a $(k-1)$-marked
part $\alpha_j$
subject to one of the following conditions:
\begin{itemize}
\item[1.] The underlying part $t$ of $\alpha_j$ satisfies
\begin{equation}\label{ft1}
f_{t}(\alpha)+f_{\overline{t}}(\alpha)+f_{t+1}(\alpha)<k-1;
  \end{equation}
\item[2.] The underlying part $t$ of $\alpha_j$ satisfies that
\begin{equation}
f_{t}(\alpha)+f_{\overline{t}}(\alpha)+f_{t+1}(\alpha)=k-1.
  \end{equation}
Moreover, we have
\begin{equation}\label{ft2}
f_{t+1}(\alpha)+f_{\overline{t+1}}(\alpha)+f_{t+2}(\alpha)<k-1.
  \end{equation}
  \end{itemize}

It is easily seen that relation (\ref{ft2}) holds for
 the largest $(k-1)$-marked part $\alpha_j$ of $\alpha$ with underlying part $t$.
This implies    there exists at least one   $(k-1)$-marked part $\alpha_j$ satisfying one of  the above two conditions.
Our goal is to find a part of $\alpha$ with underlying part  $t-1$ or $t$ and we shall increase this underlying part by one.

\noindent{\bf The Second Dilation Operation.} Let $\alpha=(\alpha_1,\ldots,\alpha_m)$ be an
overpartition in $P_{N_1,N_2,\ldots,N_{k-1};i}(n)$.
Let  $\alpha_j$ be a $(k-1)$-marked part with underlying part  $t$ for which
one of the above two conditions holds.

We first consider the case when  Condition 1 holds. Since  $t$ is the underlying part of $\alpha_j$ and
$f_{t}(\alpha)<k-1$, we deduce that
there exists a part with underlying part $t-1$. So we may assume that $r$  is the
largest mark of a part with underlying part $t-1$.
If $r=1$,  we replace the $1$-marked overlined part $\overline{t-1}$ with an $1$-marked overlined part $\overline{t}$.
If $r>1$, we replace this $r$-marked non-overline part $t-1$ with an $r$-marked part $t$.

We now consider the case when Condition 2 holds.
 In this case, we observe that there is no $k-1$-marked part with underlying part $t+1$. Moreover,
 if (\ref{ft2}) holds for $k=2$, then we replace  $\alpha_j$ with a $1$-marked part  $\overline{t+1}$.
 If (\ref{ft2}) holds for $k>2$, then
we replace  $\alpha_j$ with a $(k-1)$-marked part  $t+1$.

In either case, we obtain the  Gordon marking representation of an overpartition
 in $P_{N_1,N_2,\ldots,N_{k-1};i}(n)\setminus Q_{N_1,N_2,\ldots,N_{k-1};i}(n)$.

It can be checked that the second reduction operation  is the inverse
of the second dilation operation.
We are now ready to give a bijective proof of Theorem 5.1.

\noindent{\it Proof of Theorem \ref{PQ}.}
Using the reduction operation, we shall establish  a bijection $\psi$ between
$P_{N_1,N_2,\ldots,N_{k-1};i}$ and $Q_{N_1,N_2,\ldots,N_{k-1};i}\times R_{N_{k-1}}$,
where $R_{N_{k-1}}$ denotes the set of ordinary partitions with at most $N_{k-1}$ parts.
Let $\alpha$ be an overpartition in $ P_{N_1,N_2,\ldots,N_{k-1};i}$. Assume that
\[\alpha^{(k-1)}_1<\alpha^{(k-1)}_2<\ldots<\alpha^{(k-1)}_{N_{k-1}}\]
are the  $(k-1)$-marked parts of $\alpha$.

 Let us describe the procedure to construct $\psi(\alpha)$ by successively applying the
  second reduction operation. Keep in mind that  $\psi(\alpha)$ is
a pair $(\gamma, \delta)$, where $\gamma$ is an overpartition
in $Q_{N_1,N_2,\ldots,N_{k-1};i}$ and $\delta$ is a partition in $R_{N_{k-1}}$ such
that $|\alpha|=|\gamma|+|\delta|$.

 As discussed before, there always exists a $(k-1)$-marked part $\alpha_j$
  which satisfies either Condition 1 or Condition 2 in the second reduction operation.
  We choose the smallest $(k-1)$-marked part which satisfies either Condition 1 or
 Condition 2. Assume that it is the $l$-th $(k-1)$-marked part of $\alpha$, denoted
 $\alpha^{(k-1)}_l$. Notice that after applying  reduction operation  by choosing
 $\alpha_j$ to be $\alpha^{(k-1)}_l$, the $(l+1)$-th $(k-1)$-marked part $\alpha^{(k-1)}_{l+1}$
 remains unchanged and it satisfies the Condition 1 or Condition 2. So  can continue to apply the
  reduction operation by choosing $\alpha_j$ to be  $\alpha^{(k-1)}_{l+1}$.
Moreover, we can iterate this process with respect to the following $(k-1)$-marked
parts
$\alpha^{(k-1)}_l,\alpha^{(k-1)}_{l+1},\ldots,\alpha^{(k-1)}_{N_{k-1}}$ to get an overpartition in $Q_{N_1,N_2,\ldots,N_{k-1};i}$. Meanwhile, during the above
 process we obtain an ordinary partition with at most $N_{k-1}$ parts.

 We now give a detailed description of  the bijection $\psi$ which consists of the following steps.

\noindent Step 1. Set $\delta=\phi$ and $t=0$. We choose
the smallest $(k-1)$-marked part $\alpha^{(k-1)}_l$ which satisfies either Condition 1 or
 Condition 2. If $l=1$ and the number of parts with underlying part  $1$ is less than $i$,  go to Step 2; Otherwise, set $v=l$ and  go to Step 3.

\noindent Step 2. Recall that by the definition of $P_{N_1,N_2,\ldots,N_{k-1};i}$,
$i$ is
the maximum number of occurrences of $1$ and $\overline{1}$ in $\alpha$. There are two cases.
 If $1\leq i\leq k-1$, we repeatedly apply the reduction operation  to $\alpha$
 by choosing  $\alpha_j$ to be $\alpha^{(k-1)}_1$ until $\alpha$ becomes  an overpartition
containing an overlined part $\overline{1}$ and $i-1$ non-overlined parts $1$.
 If $i=k$, we repeatedly apply the  reduction operation to $\alpha$  by choosing  $\alpha_j$ to be $\alpha^{(k-1)}_1$ until $\alpha$ becomes an overpartition containing  an overlined part $\overline{1}$ and $k-2$ non-overlined parts $1$.  In either case, let $t$ be
the number of  the reduction operations that have been applied, and
add $t$ to $\delta$ as a new part. Set $v=2$ and go to Step 3.

\noindent Step 3. For each $s$ from $v$ to $N_{k-1}$,
 we repeatedly apply  the second reduction operations by choosing  the $(k-1)$-marked part
$\alpha_j$ to be $\alpha^{(k-1)}_s$ until $\alpha^{(k-1)}_s$ satisfies
neither Condition 1 nor Condition 2.
 After each reduction   we reset the resulting
overpartition back to $\alpha$.
Let $t$ be the number of reductions that have been applied. Add $t$ to $\delta$ as a new part.

\noindent Step 4. Let $\gamma=\alpha$ and set $\psi(\alpha)=(\gamma, \delta)$.

It can be seen that  $\gamma$ is an overpartition  in $Q_{N_1,N_2,\ldots,N_{k-1};i}$.
Meanwhile, there are $N_{k-1}-l+1$ parts in  $\delta$.
This implies that $\delta$ is a partition in $R_{N_{k-1}}$.
Moreover we have $|\alpha|=|\gamma|+|\delta|$. An example is given after the proof.

Here is an outline of the inverse of  $\psi$.
Let $\gamma$ be an overpartition in $Q_{N_1,N_2,\ldots N_{k-1};i}$ and
$\delta$
 be a  partition with $m$ parts, where $m\leq N_{k-1}$. Express the parts of  $\delta$ as
 \[\delta_1\geq \cdots\geq \delta_m.\]
 The following is a procedure to construct $\psi^{-1}(\gamma,\delta)$, which is
an overpartition $\alpha$ in $P_{N_1,N_2,\ldots N_{k-1};i}$.

\noindent Step 1. Let $\alpha=\gamma$.

\noindent Step 2.
For $t$ from $1$ to $m$,
  apply  the  dilation operation $\delta_t$ times by choosing $\alpha_j$ to be  $\gamma^{(k-1)}_{N_{k-1}-t+1}$.

\noindent Step 3.  Set $\psi^{-1}(\gamma, \delta)=\alpha$.

It can be verified that the map $\psi^{-1}(\gamma, \delta)$ is indeed the inverse
of $\psi$. The details are omitted. So we have completed the proof of Theorem \ref{PQ}. \qed

We conclude this section with  an example to demonstrate the above bijection.
For $k=4$ and $i=1$, let $\alpha$ be  an  overpartition in $P_{7,6,5;1}(128)$ as
 given by \eqref{alpha}, namely,
\begin{equation}
\setcounter{MaxMatrixCols}{13}\alpha=\begin{bmatrix}\nonumber
 \ &2&\ &\ &5&6 &\ &8&\ &\ &\ &12&\\[3pt]
\ &2&\ &4 &\ &6&\ &8&\ &10 &\ &12 &\\[3pt]
\ &\overline{2} &\ &\overline{4}&\ &\overline{6}&\overline{7}&\ &\ &\overline{10}&\overline{11} &\ &\overline{13}
\end{bmatrix}\begin{matrix}
3\\[3pt]2\\[3pt]1
\end{matrix}\;.
\end{equation}

We apply the second reduction operation by choosing $\alpha_j$  to be the $3$-marked part
 $\alpha^{(3)}_1=2$. Then $\alpha$ is mapped to an overpartition containing a part $\overline{1}$
 and no parts $1$.
Note that $i=1$. Thus we cannot further apply the
 reduction  by choosing $\alpha_j$ to be $\alpha^{(3)}_1$.
Then we get $\delta=(1)$ and $\alpha$ is  an overpartition in $P_{7,6,5;1}(127)$:

\begin{equation}
\setcounter{MaxMatrixCols}{13}\alpha=\begin{bmatrix}\nonumber
 \ &2&\ &\ &5&6 &\ &8&\ &\ &\ &12&\\[3pt]
\ &2&\ &4 &\ &6&\ &8&\ &10 &\ &12 &\\[3pt]
\overline{1}&\ &\ &\overline{4}&\ &\overline{6}&\overline{7}&\ &\ &\overline{10}&\overline{11} &\ &\overline{13}
\end{bmatrix}\begin{matrix}
3\\[3pt]2\\[3pt]1
\end{matrix}\;.
\end{equation}

Next we choose $\alpha_j$ to be  $\alpha^{(3)}_2$. Then we can apply  reduction three times to
change the $3$-marked part $5$ to  the $3$-marked part $4$, change  the $1$-marked part $\overline{4}$ to the $1$-marked part $\overline{3}$, and change
the $2$-marked part $4$ to the $2$-marked part $3$. After that  $\alpha^{(3)}_2$ no longer
satisfies   Condition 1 or Condition 2.  Then we add $3$ to $\delta$ as a new part to get $\delta=(3,1)$ and
 $\alpha$ becomes an overpartition in $P_{7,6,5;1}(124)$:
\begin{equation}
\setcounter{MaxMatrixCols}{13}\begin{bmatrix}\nonumber
 \ &2&\ &4 &\ &6 &\ &8&\ &\ &\ &12&\\[3pt]
\ &2&3 &\ &\ &6&\ &8&\ &10 &\ &12 &\\[3pt]
\overline{1}&\ &\overline{3} &\ &\ &\overline{6}&\overline{7}&\ &\ &\overline{10}&\overline{11} &\ &\overline{13}
\end{bmatrix}\begin{matrix}
3\\[3pt]2\\[3pt]1
\end{matrix}\;.
\end{equation}

We continue to consider $\alpha^{(3)}_3=6$ as a choice of $\alpha_j$.
We can apply   reduction three times so that $\alpha$ becomes an overpartition in $P_{7,6,5;1}(121)$  as given below:
\begin{equation}
\setcounter{MaxMatrixCols}{13}\begin{bmatrix}\nonumber
 \ &2&\ &4 &5 &\ &\ &8&\ &\ &\ &12&\\[3pt]
\ &2&3 &\ &5 &\ &\ &8&\ &10 &\ &12 &\\[3pt]
\overline{1}&\ &\overline{3} &\ &\overline{5} &\ &\overline{7}&\ &\ &\overline{10}&\overline{11} &\ &\overline{13}
\end{bmatrix}\begin{matrix}
3\\[3pt]2\\[3pt]1
\end{matrix}\;.
\end{equation}
Then add $3$ as a new part to $\delta$ and get $\delta=(3,3,1)$.

For the remaining  $3$-marked parts $8$
we can apply the reduction three times by choosing $\alpha_j=8$. Finally,
 for the $3$-marked part $12$, we can apply the reduction seven times
by choosing $\alpha_j=12$. Thus  we get
 $\delta=(7,3,3,3,1)$. In the mean time, $\alpha$  is mapped to an  overpartition in $Q_{7,6,5;1}(111)$ as given by
\begin{equation}
\gamma=\setcounter{MaxMatrixCols}{13}\begin{bmatrix}\nonumber
 \ &2&\ &4 &5 &\ &7 &\ &9 &\ &\ &&\\[3pt]
\ &2&3 &\ &5 &\ &7 &8&\ &\ & &12 &\\[3pt]
\overline{1}&\ &\overline{3} &\ &\overline{5} &\overline{6} &\ &\overline{8} &\ &\ &\overline{11} &\ &\overline{13}
\end{bmatrix}\begin{matrix}
3\\[3pt]2\\[3pt]1
\end{matrix}\; .
\end{equation}

\section{The third bijection for the proof of Theorem \ref{main}}

In this section,  we   give the third bijection for the proof of Theorem \ref{main},
which is  between  $Q_{N_1,\ldots,N_{k-1};i}$ and $Q_{N_1-1,\ldots,N_{k-1}-1;i}$. By
this correspondence, we can derive
  a recurrence relation on $Q_{N_1,\ldots,N_{k-1};i}$,
   which yields the generating function of $Q_{N_1,\ldots,N_{k-1};i}$ as stated in
   the following theorem.

\begin{thm}\label{gfQ}
For $k\geq 2$ and $1\leq i\leq k$, we have
\begin{equation}\label{eqQ}
 \sum_{\gamma \in Q_{N_1,\ldots,N_{k-1};i}}x^{l(\gamma)}q^{|\gamma|}=
\frac{q^{\frac{(N_1+1)N_1}{2}+N_2^2+\cdots+N_{k-1}^2+N_{i+1}+\cdots+N_{k-1}}x^{N_1+\cdots+N_{k-1}}}
{(q)_{N_1-N_2}
\cdots(q)_{N_{k-2}-N_{k-1}}}.
\end{equation}
\end{thm}

In order to prove the above theorem by induction, we need the following bijection.

\begin{thm}\label{indQ} For $N_{k-1}>0$,
there is  a bijection between $Q_{N_1,\ldots,N_{k-1};i}(n)$ and
$Q_{N_1-1,\ldots,N_{k-1}-1;i}(n-N_1-2N_2-\cdots-2N_{k-1}+i-1)$. In terms of generating
functions, we have
\begin{equation}
 \sum_{\gamma \in Q_{N_1,\ldots,N_{k-1};i}}q^{|\gamma|}=
 q^{N_1+2N_2+\ldots+2N_{k-1}-i+1}
\sum_{\gamma \in Q_{N_1-1,\ldots,N_{k-1}-1;i}}q^{|\gamma|}.
\end{equation}
\end{thm}

\pf Assume that $N_{k-1}>0$. We will give a bijection $\chi$ between $Q_{N_1,\ldots,N_{k-1};i}(n)$ with
and $Q_{N_1-1,\ldots,N_{k-1}-1;i}(n-N_1-2N_2-\cdots-2N_{k-1}+i-1)$.
Let $\gamma$ be an overpartition in $Q_{N_1,\ldots,N_{k-1};i}(n)$.
We proceed to construct $\chi(\gamma)$, which is
an overpartition  $\mu$ in $Q_{N_1-1,\ldots,N_{k-1}-1;i}(n-N_1-2N_2-\cdots-2N_{k-1}+i-1)$.

The idea of this bijection goes as follows.
For each $1$-marked part $\gamma_j$ with underlying part $a_j$,
we shall allocate  a part with underlying part $a_j$ subject to certain conditions.
Then we  increase this part by $1$. Furthermore, for
 each $1$-marked part, we remove  the smallest part of  each row in the Gordon marking representation of  the resulting overpartition, and subtract $2$ from the other parts.
Here are the detailed description.

\noindent Step 1. Let $\mu=\gamma$.

\noindent Step 2.
For $i$ from $N_1$ to $1$,
let $t$ be the underlying part of  $\mu^{(1)}_i$.

If
there are two  parts  of the same mark but with distinct underlying parts   $t-1$ and $t$,
we denote this mark by $r$. Then we change the
$r$-marked part  with underlying part  $t$ to an $r$-marked part with underlying part $t+1$;

Otherwise, we find the greatest  mark $r$, such that there is an $r$-marked part with underlying part
$t$.  If $r=1$, replace the $1$-marked overlined part  $\overline{t}$  of $\mu$
 with an $1$-marked part $\overline{t+1}$. If $r>1$, replace the
$r$-marked part $t$  with an $r$-marked part $t+1$.
Clearly, the sum of the parts of $\mu$ becomes $n+N_1$.

\noindent Step 3. Delete $\mu^{(1)}_1,\ldots,\mu^{(k-1)}_1$ and
subtract $2$ from each part of $\mu$.

From the definition of $Q_{N-1,\ldots,N_{k-1};i}$, the smallest part of each row is $1$ or $2$.
Clearly, after Step 2 there are $i-1$ parts equal to $1$ and $k-i$ parts equal to $2$ in $\mu$. So after Step 3 the sum of parts of $\mu$ equals
\[n+N_1-(i-1)-2(k-i)-2(N_1+\cdots+N_{k-1}-(k-1))=n-N_1-2N_2-\cdots-2N_{k-1}+i-1.\]

\noindent Step 4. Let $\chi(\lambda)=\mu$.

It can be seen that after the above process  we obtain the Gordon marking of an overpartition in
$Q_{N_1-1,\ldots,N_{k-1}-1;i}(n-N_1-2N_2-\cdots-2N_{k-1}+i-1)$.

We now consider the inverse of  $\chi$.
Let $\mu \in Q_{N_1-1,\ldots,N_{k-1}-1;i}(n)$.
The following is a procedure to construct $\chi^{-1}(\mu)$, which is
a partition $\gamma$ in $Q_{N_1-1,\ldots,N_{k-1}-1;i}(n+N_1+2N_2+\ldots+2N_{k-1}-i+1)$.

\noindent Step 1. Let $\gamma=\mu$.

\noindent Step 2. Increase each part of $\gamma$ by $2$.

\noindent Step 3.
If $i=1$, we add $1$-marked part $\overline{2}$, a $2$-marked part $2$, $\ldots$, and a $(k-1)$-marked part $2$  to
$\gamma$ as new parts.
If $i\geq 2$, we add a $1$-marked part $\overline{1}$, $\ldots$, an $(i-1)$-marked part $1$,
an $i$-marked part $2$, $\ldots$,  and a  $(k-1)$-marked part $2$  to $\gamma$  as new parts.
Now $\gamma$ contains $N_{1}+1$ parts with $1$-marked.

\noindent Step 4. For $j$ from 1 to $N_{1}+1$, let $t$ be the underlying part of
 $\gamma^{(1)}_j$.

If  $\overline{t+1}$  is a part of $\gamma$ or there are no parts with underlying part $t+1$,
then we replace the  overlined $1$-marked  part $\overline{t}$ with a $1$-marked part $\overline{t-1}$;

If
$\overline{t+1}$ is not a part of $\gamma$ but $t+1$ is a part of $\gamma$, then we choose the
smallest mark $r$ of parts with underlying part $t+1$, and replace this $r$-marked part $t+1$ with an $r$-marked part $t$.

\noindent Step 5. Set $\chi^{-1}(\mu)=\gamma$.

It can be verified that after the above steps we get the Gordon marking
of an overpartition in
$Q_{N_1-1,\ldots,N_{k-1}-1;i}(n+N_1+2N_2+\ldots+2N_{k-1}-i+1)$.

It is routine to check that the map $\chi^{-1}$ is the inverse of $\chi$.
\qed

Here we give an example of the above  bijection.
Let
$\gamma=\\
(\overline{1},1,\overline{2},2,2,\overline{3},3,3,4,4,\overline{5},5,5,5,\overline{6},6,7,7,7,\overline{8},
8,8,\overline{9},9,9,10,10,\overline{11}, 11,12,12,\overline{14},14,15,\overline{17},17,17)$
in $Q_{10,9,8,6,6;1}(311)$. Set $\mu=\gamma$.
Below is the Gordon marking representation of $\mu$
\begin{equation}
\setcounter{MaxMatrixCols}{17}\begin{bmatrix}
\ &2&\ &4&\ &6&\ &8&\ &10&\ &12&\ &\ &\ &\ &\\[3pt]
\ &2&\ &4&\bf{5}&\ &7&\ &9&\ &11&\ &\ &\ &\ &\ &\\[3pt]
\ &2&\bf{3}&\ &5&\ &7&\ &9&\ &11&\ &\ &\ &15&\ &\bf{17}\\[3pt]
\bf{1}&\ &3&\ &5&\ &7&\bf{8}&\ &10&\bf{11}&\ &\ &14&\ &\ &17\\[3pt]
\overline{1}&\ &\overline{3}&\ &\overline{5}&\bf{\overline{6}}&\ &\overline{8}&\bf{\overline{9}}&\ &\overline{11}&\
&\bf{\overline{13}}&\bf{\overline{14}}&\ &\ &\overline{17}
\end{bmatrix}
\begin{matrix}
5\\[3pt]4\\[3pt]3\\[3pt]2\\[3pt]1
\end{matrix}\; ,
\end{equation}
where the parts in boldface are those we should move to the right in Step 2.
After Step 2, $\mu$ is changed to
\begin{equation}
\setcounter{MaxMatrixCols}{20}\begin{bmatrix}\nonumber
 \ &2&\ &4&\ &6&\ &8&\ &10&\ &12&\ &\ &\ &\ &\ &\ \\[3pt]
 \ &2&\ &4 &\ &\mathbf{6}&7 &\ &9 &\ &11&\ &\ &\ &\ &\ &\ &\ \\[3pt]
 \ &2&\ &\mathbf{4}&5&\ &7&\ &9&\ &11&\ &\ &\ &15&\ &\ &\mathbf{18}\\[3pt]
\ &\mathbf{2}&3 &\ &5 &\ &7 &\ &\mathbf{9}&10&\ &\mathbf{12}&\ &14&\ &\ &17&\ \\[3pt]
\overline{1}&\ &\overline{3} &\ &\overline{5} &\  &\mathbf{\overline{7}} &\overline{8} &\
&\mathbf{\overline{10}} &\overline{11} &\ &\ &\mathbf{\overline{14}}&\mathbf{\overline{15}}&\ &\overline{17}
\end{bmatrix}\begin{matrix}5\\[3pt]4\\[3pt]3\\[3pt]2\\[3pt]1\end{matrix}\; .
\end{equation}
Deleting the parts $\mu^{(1)}_1,\ldots,\mu^{(5)}_1$ and subtracting $2$ from the
other parts of $\mu$, we get
\begin{equation}
\setcounter{MaxMatrixCols}{16}\begin{bmatrix}
\ &2&\ &4&\ &6&\ &8&\ &10&\ &\ &\ &\ &\ &\\[3pt]
\ &2&\ &4&5&\ &7&\ &9&\ &\ &\ &\ &\ &\ &\\[3pt]
\ &2&3&\ &5&\ &7&\ &9&\ &\ &\ &13&\ &\ &16\\[3pt]
1&\ &3&\ &5&\ &7&8&\ &10&\ &12 &\ &\ &15 &\\[3pt]
\overline{1}&\ &\overline{3} &\ &\overline{5}&\overline{6}&\ &\overline{8}&\overline{9}&\ &\ &\overline{12}
&13&\ &\overline{15} &
\end{bmatrix}
\begin{matrix}
5\\[3pt]4\\[3pt]3\\[3pt]2\\[3pt]1
\end{matrix}\; ,
\end{equation}
which is the Gordon marking representation of an overpartition in $Q_{9,8,7,5,5;1}(254)$.
It can be checked that the above process is reversible.

{\noindent\it The proof of Theorem \ref{gfQ}.} We use induction on $k$.
For $k=2$ and $i=1$, the generating function of $Q_{N_1;1}$ is
\[\sum_{\lambda \in Q_{N_1;1}}q^{|\lambda|}=
q^{\frac{(N_1+1)N_1}{2}}.\]
For $k=2$ and $i=2$, the generating function of  $Q_{N_1;2}$ is
\[\sum_{\lambda \in Q_{N_1;2}}q^{|\lambda|}=
q^{\frac{(N_1+1)N_1}{2}}.\]
So Theorem \ref{gfQ} holds for $k=2$. Assume that it holds for $k-1$,
   that is,
\[
 \sum_{\lambda \in Q_{N_1,\ldots,N_{k-2};i}}q^{|\lambda|}=
\frac{q^{\frac{(N_1+1)N_1}{2}+N_2^2+\cdots+N_{k-2}^2+N_{i+1}+\cdots+N_{k-2}}}{(q)_{N_1-N_2}(q)_{N_2-N_3}
\cdots(q)_{N_{k-3}-N_{k-2}}}.
\]
We proceed to show that it holds for $Q_{N_1,\ldots,N_{k-1};i}$.

If $N_{k-1}=0$, by the definitions of $Q_{N_1,\ldots,N_{k-2},0;i}$ and $P_{N_1,\ldots,N_{k-2};i}$, we find that
\[Q_{N_1,\ldots,N_{k-2},0;i}=P_{N_1,\ldots,N_{k-2};i}.\]
In view of Theorem \ref{PQ}, the generating function of $Q_{N_1,\ldots,N_{k-2},0;i}$ equals
\begin{equation}\label{gfQ0}\sum_{\lambda \in Q_{N_1,\ldots,N_{k-2},0;i}}q^{|\lambda|}=
\frac{1}{(q)_{N_{k-2}}}\times\frac{q^{\frac{(N_1+1)N_1}{2}+N_2^2
+\cdots+N_{k-2}^2+N_{i+1}+\cdots+N_{k-2}}}{(q)_{N_1-N_2}(q)_{N_2-N_3}
\cdots(q)_{N_{k-3}-N_{k-2}}}.\end{equation}

If $N_{k-1}>0$,  applying  Theorem \ref{indQ} $N_{k-1}$ times, we obtain that
\begin{align} \label{2n1}
 &\sum_{\lambda \in Q_{N_1,\ldots,N_{k-1};i}}q^{|\lambda|}\nonumber\\
&\quad =q^{\frac{(2N_1-N_{k-1}+1)N_{k-1}}{2}+(2N_2-N_{k-1}+1)N_{k-1}+\cdots+(N_{k-1}+1)N_{k-1}-N_{k-1}i+N_{k-1}}\nonumber\\
&\qquad  \quad  \times\sum_{\lambda \in Q_{N_1-N_{k-1},\ldots,N_{k-2}-N_{k-1},0;i}}q^{|\lambda|}.
 \end{align}
Combining  \eqref{gfQ0} and  (\ref{2n1}), we have for $1\leq i \leq k-1$
\begin{align*}
&\sum_{\lambda \in Q_{N_1,\ldots,N_{k-1};i}}q^{|\lambda|}\nonumber\\[6pt]
&\quad =q^{\frac{(2N_1-N_{k-1}+1)N_{k-1}}{2}+(2N_2-N_{k-1}+1)N_{k-1}+\cdots+(N_{k-1}+1)N_{k-1}-N_{k-1}i+N_{k-1}}\\
&\qquad \quad \times\sum_{\lambda \in Q_{N_1-N_{k-1},\ldots,N_{k-2}-N_{k-1},0;i}}q^{|\lambda|}\nonumber\\[6pt]
&\quad =\frac{q^{\frac{(N_1+1)N_1}{2}+N_2^2+\cdots+N_{k-1}^2+N_{i+1}+\cdots+N_{k-1}}}{(q)_{N_1-N_2}(q)_{N_2-N_3}
\cdots(q)_{N_{k-2}-N_{k-1}}}.
 \end{align*}

Since for any overpartition in $Q_{N_1,\ldots,N_{k-1};i}$ the  smallest $1$-marked part is overlined,  the non-overlined $1$ can occur at most $k-2$ times. This implies that
$Q_{N_1,\ldots,N_{k-1};k}=Q_{N_1,\ldots,N_{k-1};k-1}$.
we have proved that identity \eqref{eqQ} holds for $1\leq i\leq k$, that is,
Theorem \ref{gfQ} holds for $k$. This completes the proof.
\qed

\section{Proof of Theorem \ref{main}}

In this section, we finish the proof of Theorem \ref{main}.
Using the three bijections given in the previous sections, we
can derive the generating function  of $F_{k,i}(m,n)$ as stated in Theorem
\ref{F}.  Then we   compute the generating function of $G_{k,i}(m,n)$
which leads to the generating function of $D_{k,i}(m,n)$.
We first give the proof of Theorem \ref{F}.

\noindent
{\it Proof of Theorem \ref{F}.} By  Theorems \ref{FP}, \ref{PQ}, and \ref{gfQ}, we find that  the generating function of $F_{k,i}(m,n)$ equals
\begin{align*}
&\sum_{n=0}^{\infty}F_{k,i}(m,n)x^mq^n\\&\qquad=\sum_{N_1\geq N_2\geq\cdots\geq
N_{k-1}\geq0}\frac{(-q)_{N_1-1}}{(q)_{N_{k-1}}}\sum_{\lambda \in Q_{N_1,\ldots,N_{k-1};i}}x^{N_1+\cdots+N_{k-1}}q^{|\lambda|}
\\[6pt] &\qquad=\sum_{N_1\geq N_2\geq\cdots\geq
N_{k-1}\geq0}\frac{q^{\frac{(N_1+1)N_1}{2}+N_2^2+\cdots+N_{k-1}^2+N_{i+1}+\cdots+N_{k-1}}
(-q)_{N_1-1}x^{N_1+\cdots+N_{k-1}}}{(q)_{N_1-N_2}\cdots(q)_{N_{k-2}-N_{k-1}}(q)_{N_{k-1}}},\end{align*}
as claimed. \qed

 Given the relation between $F_{k,i}(m,n)$ and $G_{k,i}(m,n)$ as stated in Lemma \ref{FG1}, we can derive the generating function of $G_{k,i}(m,n)$.

\begin{thm}\label{G}
For $k\geq i\geq 1$,
\begin{align}\label{G0}
&\sum_{n=0}^{\infty}G_{k,i}(m,n)x^mq^n\nonumber
\\&\quad=\sum_{N_1\geq N_2\geq\cdots\geq
N_{k-1}\geq0}\frac{q^{\frac{(N_1+1)N_1}{2}+N_2^2+\cdots+N_{k-1}^2+N_i+\cdots+N_{k-1}}
(-q)_{N_1-1}x^{N_1+\cdots+N_{k-1}}}{(q)_{N_1-N_2}\cdots(q)_{N_{k-2}-N_{k-1}}(q)_{N_{k-1}}}.\end{align}
 \end{thm}

\pf  From relation \eqref{FG}, we deduce that for $2\leq i\leq k$,
\begin{align}\label{G1}&\sum_{m,n\geq 0}G_{k,i}(m,n)x^mq^n\nonumber\\&\qquad=\sum_{m,n\geq 0}F_{k,i-1}(m,n)x^mq^n\nonumber\\&
\qquad =\sum_{N_1\geq N_2\geq\cdots\geq
N_{k-1}\geq0}\frac{q^{\frac{(N_1+1)N_1}{2}+N_2^2+\cdots+N_{k-1}^2+N_i+\cdots+N_{k-1}}
(-q)_{N_1-1}x^{N_1+\cdots+N_{k-1}}}{(q)_{N_1-N_2}\cdots(q)_{N_{k-2}-N_{k-1}}(q)_{N_{k-1}}}.
\end{align}
For $i=1$, from  \eqref{FG2} it follows that
 \[\sum_{m,n\geq 0}G_{k,1}(m,n)x^mq^n=\sum_{m,n \geq 0}F_{k,k}(m,n)(xq)^mq^n.\]
Using the generating function of $F_{k,k}(m,n)$, we obtain
\begin{align}\label{G2}&\sum_{m,n\geq 0}G_{k,1}(m,n)x^mq^n\nonumber
\\&\qquad =\sum_{N_1\geq N_2\geq\cdots\geq
N_{k-1}\geq0}\frac{q^{\frac{(N_1+1)N_1}{2}+N_2^2+\cdots+N_{k-1}^2+N_1+\cdots+N_{k-1}}
(-q)_{N_1-1}x^{N_1+\cdots+N_{k-1}}}{(q)_{N_1-N_2}\cdots(q)_{N_{k-2}-N_{k-1}}(q)_{N_{k-1}}}.\end{align}
Observe that the above formulas \eqref{G1} for $i>1$ and \eqref{G2} for $i=1$ take the same form \eqref{G0} as
in the theorem. This completes the proof.
 \qed

We are now ready to finish the proof of Theorem \ref{main}.

\noindent
{\it Proof of Theorem \ref{main}.} By the generating functions of $G_{k,i}(m,n)$ and $F_{k,i}(m,n)$ and relation \eqref{DFG}, we find that
\begin{align*}&\sum_{m,n\geq 0}D_{k,i}(m,n)x^mq^n\\&\quad=\sum_{m,n\geq 0}F_{k,i}(m,n)x^mq^n+\sum_{m,n\geq 0}G_{k,i}(m,n)x^mq^n\\&\quad=\sum_{N_1\geq N_2\geq\cdots\geq
N_{k-1}\geq0}\frac{q^{\frac{(N_1+1)N_1}{2}+N_2^2+\cdots+N_{k-1}^2+N_{i+1}+\cdots+N_{k-1}}
(-q)_{N_1-1}(1+q^{N_i})x^{N_1+\cdots+N_{k-1}}}{(q)_{N_1-N_2}\cdots(q)_{N_{k-2}-N_{k-1}}(q)_{N_{k-1}}}.
\end{align*}
This completes the proof of Theorem \ref{main}. \qed

\vspace{0.5cm}
 \noindent{\bf Acknowledgments.}  This work was supported by  the 973
Project, the PCSIRT Project of the Ministry of Education,  and the
National Science Foundation of China.

\end{document}